\documentclass[leqno,12pt]{article}
\usepackage{amsmath,amssymb,wasysym,rotating,a4wide,color}

\usepackage[all,cmtip]{xy}

\usepackage{verbatim}

\newdir^{ (}{{}*!/-3pt/\dir^{(}}
\newdir_{ (}{{}*!/-3pt/\dir_{(}}

\entrymodifiers={+!!<0pt,\fontdimen22\textfont2>}

\def\sgn{\mathop{\rm sgn}\nolimits}
\def\bigtimes{\mathop{\raise-2pt\hbox{\huge$\times$}}}

\newbox\circbulletbox
\setbox\circbulletbox=\hbox{\,{\large$\circledcirc$}\kern-8.7pt\raise .6pt\hbox{$\bullet$}\,}
\def\circbullet{\copy\circbulletbox}
\def\myothersmash#1#2{\smash{\raisebox{#1}{$#2$}}}

\let\le\leqslant
\let\ge\geqslant

\def\mycirc{{\kern1pt\circ\kern2pt}}

\def\Aut{\mathop{\rm Aut}\nolimits}

\def\Gal{\mathop{\rm Gal}\nolimits}

\def\Spec{\mathop{\rm Spec}\nolimits}

\def\Ker{\mathop{\rm Ker}\nolimits}

\def\Norm{\mathop{\rm Norm}\nolimits}

\def\diag{\mathop{\rm diag}\nolimits}

\def\proj{\mathop{\rm pr}\nolimits}
\def\GL{\mathop{\rm GL}\nolimits}

\def\et{{\rm \acute{e}t}}

\def\geom{{\rm geom}}
\def\arith{{\rm arith}}

\let\phi\varphi
\let\theta\vartheta
\let\epsilon\varepsilon
\let\setminus\smallsetminus

\newtheorem{Thm}{Theorem}[section]
\newtheorem{Prop}[Thm]{Proposition}
\newtheorem{Lem}[Thm]{Lemma}

\newtheorem{Rem}[Thm]{Remark}
\newtheorem{Ex}[Thm]{Example}

\numberwithin{Thm}{section}
\def\UseTheoremCounterForNextEquation{\setcounter{equation}{\value{Thm}}\addtocounter{Thm}{1}}
\renewcommand{\theequation}
             {\arabic{section}.\arabic{Thm}}

\def\qed{{\hskip0pt\unskip\unskip\nobreak\hfil\penalty50
          \hskip1em\hbox{}\nobreak\hfil
           {$\square$}
          \parfillskip=0pt\finalhyphendemerits=0
          \par}\medskip}

\newenvironment{Proof}
               {\noindent{\bf Proof.}\ }
               {\qed}

               {\noindent{\bf Proof of #1.}\ }
               {\qed}

\newcommand{\BA}{{\mathbb{A}}}

\newcommand{\BF}{{\mathbb{F}}}

\newcommand{\BP}{{\mathbb{P}}}
\newcommand{\BQ}{{\mathbb{Q}}}

\newbox\mybox
\def\arrover#1{\mathrel{
       \setbox\mybox=\hbox spread 1.4em
              {\hfil$\scriptstyle#1$\hfil}
       \vbox{\offinterlineskip\copy\mybox
             \hbox to\wd\mybox{\rightarrowfill}}}}

\def\larrover#1{\mathrel{
       \setbox\mybox=\hbox spread 1.4em
              {\hfil$\scriptstyle#1\vphantom{g}$\hfil}
       \vbox{\offinterlineskip\copy\mybox
             \hbox to\wd\mybox{\leftarrowfill}}}}

\def\ontoover#1{\mathrel{
       \setbox\mybox=\hbox spread 1.4em
              {\hfil$\scriptstyle#1\vphantom{g}$\hfil}
       \vbox{\offinterlineskip\copy\mybox
             \hbox to\wd\mybox{\rightarrowfill\hskip-2.8mm
                               $\rightarrow$}}}}
\def\leftontoover#1{\mathrel{
       \setbox\mybox=\hbox spread 1.4em
              {\hfil$\scriptstyle#1\vphantom{g}$\hfil}
       \vbox{\offinterlineskip\copy\mybox
             \hbox to\wd\mybox{$\leftarrow$\hskip-2.8mm
                               \leftarrowfill}}}}
\let\longto\longrightarrow
\let\into\hookrightarrow
\let\onto\twoheadrightarrow

\begin{document}

\title{Profinite iterated monodromy groups\\
arising from quadratic morphisms\\
with infinite postcritical orbits}

\author{Richard Pink\\[12pt]
\small Department of Mathematics \\[-3pt]
\small ETH Z\"urich\\[-3pt]
\small 8092 Z\"urich\\[-3pt]
\small Switzerland \\[-3pt]
\small pink@math.ethz.ch\\[12pt]}


\date{September 23, 2013}

\maketitle

\begin{abstract}
We study in detail the profinite group $G$ arising as geometric \'etale iterated monodromy group of an arbitrary quadratic morphism $f$ with an infinite postcritical orbit over a field of characteristic different from two. This is a self-similar closed subgroup of the group of automorphisms of a regular rooted binary tree. In many cases it is equal to the automorphism group of the tree, but there remain some interesting cases where it is not. In these cases we prove that the conjugacy class of $G$ depends only on the combinatorial type of the postcritical orbit of~$f$. We also determine the Hausdorff dimension and the normalizer of~$G$. This result is then used to describe the arithmetic \'etale iterated monodromy group of~$f$. 

The methods used mostly group theoretical and of the same type as in a previous article of the same author dealing with quadratic polynomials with a finite postcritical orbit. The results on abstract self-similar profinite groups acting on a regular rooted binary tree may be of independent interest.
\end{abstract}

{\renewcommand{\thefootnote}{}
}

\newpage
\include{FiniPart1}
\tableofcontents


\addtocounter{section}{-1}
\section{Introduction}
\label{Intro}

This article is a sequel to the article \cite{Pink2013b} by the same author. The situation is essentially the same as there, except that here we consider a quadratic morphism which is not necessarily a polynomial and assume that its postcritical orbit is infinite. We assume familiarity with Section 1 of \cite{Pink2013b}, using the same setup and notation as there. For easier reference we envelop all the material of the present article in one Section \ref{4Infinite}.

\medskip
So let $f$ be a rational function of degree two in one variable, with coefficients in a field $k$ of characteristic different from~$2$. Let $\bar k$ be a separable closure of~$k$, and let $C\subset \BP^1(\bar k)$ denote the set of critical points of~$f$. Let $P := \bigcup_{n\ge1}f^n(C) \subset \BP^1(\bar k)$ denote the (strictly) \emph{postcritical orbit} of~$f$, which we now assume to be infinite.
Let $T$ be a regular rooted binary tree, and let $\rho\colon \pi_1^\et(\BP^1_k \setminus P) \to\Aut(T)$ be the monodromy representation which describes the infinite tower of coverings $\ldots \smash{\stackrel{f}{\longto}}\BP^1_k \smash{\stackrel{f}{\longto}}\BP^1_k \smash{\stackrel{f}{\longto}}\BP^1_k$. Its image is called the \emph{arithmetic iterated monodromy group} associated to~$f$. We are also interested in the image $G^\geom$ of the subgroup $\pi_1^\et(\BP^1_{\bar k} \setminus P)$, which is called the \emph{geometric iterated monodromy group} associated to~$f$.
As in \cite[\S1.7]{Pink2013b} the group $G^\geom$ is a pro-$2$-group which is topologically generated by elements $b_p$ for all $p\in P$, where $b_p$ is a generator of the image of some inertia group above~$p$. Moreover, a product of the generators in some order converges to the identity element in $\Aut(T)$, and each generator $b_p$ is conjugate under $\Aut(T)$ to
\UseTheoremCounterForNextEquation
\begin{equation}\label{0GenRecConj1}
\qquad\left\{\begin{array}{ll}
\ \ \ \sigma & \hbox{if $p=f(c)$ for some $c\in C\setminus P$,}\\[3pt]
(b_c,1)\,\sigma & \hbox{if $p=f(c)$ for some $c\in C\cap P$,}\\[3pt]
(b_q,1) & \hbox{if $p=f(q)$ for a unique $q\in P\setminus C$,} \\[3pt]
(b_q,b_{q'}) & \hbox{if $p=f(q)=f(q')$ for distinct $q,q'\in P\setminus C$,}
\end{array}\right\}
\end{equation}
see \cite[Prop.\ \ref{17GenRec1}]{Pink2013b} and (\ref{4GenRecConj1}). Again, all the results about $G^\geom$ in this article are purely algebraic consequences of these facts.
Once we know $G^\geom$, we describe its normalizer $N$ and determine $G^\arith$ from the natural homomorphism $\bar\rho\colon \Gal(\bar k/k)\onto G^\arith/G^\geom \subset N/G^\geom$ induced by~$\rho$.

\medskip
Let $p_0$ and $q_0$ denote the two critical points of~$f$. Our main results are:
\begin{enumerate}
\item[$\bullet$] If $f^{r+1}(p_0)\not=f^{r+1}(q_0)$ for all $r\ge1$, then $G^\geom=G^\arith=\Aut(T)$: see Theorem \ref{48GgeomThm} (a).

\item[$\bullet$] Otherwise let $r\ge1$ be minimal with $f^{r+1}(p_0)=f^{r+1}(q_0)$. Then $G^\geom$ is conjugate under $\Aut(T)$ to a certain closed subgroup $G(r)$ of $\Aut(T)$ that depends only on~$r$ and is defined by explicit recursively defined generators: see Theorem \ref{48GgeomThm} (b).

\item[$\bullet$] In either case the subgroup $G^\geom\subset\Aut(T)$ up to conjugacy depends only on the combinatorial type of the postcritical orbit of~$f$.

\item[$\bullet$] In the second case we have $G^\arith=G^\geom$ if $p_0$ and $q_0$ are defined over~$k$, otherwise $G^\arith$ is a certain extension of index~$2$ of $G^\geom$ that we describe by an additional explicit generator: see Corollary \ref{48GarithThm2}.

\item[$\bullet$] 
As a variant consider a finite extension $k'$ of $k$ and a point $x'\in\BP^1(k')\setminus P$. Let $G_{x'}$ denote the image of the Galois representation $\Gal(\bar k/k')\to \Aut(T)$ describing the action on all preimages $\coprod_{n\ge0}f^{-n}(x')$. Then by general principles there exists an inclusion $G_{x'}\subset wG^\arith w^{-1}$ for some $w\in\Aut(T)$, which in the second case above yields a nontrivial upper bound for~$G_{x'}$.
\end{enumerate}

The greater part of the article is a study of the abstractly defined group $G(r)\subset \Aut(T)$ mentioned above, which besides its direct consequences for $G^\geom$ is interesting in its own right. The methods used for this are standard finite and profinite group theory, and time and again the exploitation of the self-similarity properties of~$G(r)$ resulting from the recursion relations of its generators. 
Among other things: 
\begin{enumerate}
\item[$\bullet$] We show that the Hausdorff dimension of $G(r)$ is $1-2^{-r}$: see Theorem \ref{43Hausdorff}.

\item[$\bullet$] We prove a \emph{semirigidity} property of the generators~$G(r)$. Actually we show that for any collection of elements of $\Aut(T)$ satis\-fy\-ing certain weak recursion relations like those in (\ref{0GenRecConj1}) and whose infinite product is~$1$, on conjugating them by the same element of $\Aut(T)$ we can make them all lie in~$G(r)$, be conjugate to the standard generators of $G(r)$ under~$G(r)$, and topologically generate~$G(r)$: see Theorem \ref{44SemiRigid}. This semirigidity is the key towards identifying $G^\geom$.

\item[$\bullet$] We determine the normalizer $N(r)\subset\Aut(T)$ of $G(r)$ and describe it using further explicit generators. We construct an isomorphism between $N(r)/G(r)$ and a countably infinite product $\BF_2^{\raise2pt\hbox{$\kern1pt\scriptstyle\infty$}}$ of copies of the cyclic group of order~$2$. 

\item[$\bullet$] We determine all possible inclusions between the groups $G(r)$ and $N(r)$ for different values of~$r$: see Subsection \ref{47Inclusions}.
\end{enumerate}

\medskip
For relations with the existing literature, without any attempt at completeness, see, besides the references in \cite{Pink2013b}, the paper by Jones-Manes \cite{JonesManes} and the survey paper Jones \cite{JonesSurvey2013}.


\numberwithin{Thm}{subsection}
\def\UseTheoremCounterForNextEquation{\setcounter{equation}{\value{Thm}}\addtocounter{Thm}{1}}
\renewcommand{\theequation}
             {\arabic{section}.\arabic{subsection}.\arabic{Thm}}


\setcounter{section}{3}
%
%

\section{Infinite case}
\label{4Infinite}

Throughout this article we use the same notation and conventions as in Section 1 of \cite{Pink2013b}. 
In particular $T$ is a fixed regular rooted binary tree, and its automorphism group is denoted~$W$. The automorphism group of the truncation $T_n$ of $T$ at level $n$ is denoted~$W_n$. The symbol $\sim$ always means conjugacy under~$W$, while conjugacy under subgroups of $W$ will always be expressed in words. We heavily rely on the fundamental construction principle for elements of $W$ by recursion relations and the properties thereof explained in \cite[\S1.4]{Pink2013b}. There are natural sign homomorphisms $\sgn_n\colon W\to \{\pm1\}$ for all $n\ge1$: see \cite[\S1.5]{Pink2013b}. We also use the description of $G^\geom$ and $G^\arith$ and of generators of the former from \cite[\S1.7]{Pink2013b}.

\medskip
Subsections \ref{41TheGroups} through \ref{47Inclusions} deal with purely combinatorially defined closed subgroups of $W$ whose generators satisy weak recursion relations motivated by the shape of an infinite postcritical orbit of a quadratic morphism. Subsection \ref{41TheGroups} covers the cases where the group turns out to be~$W$.  The other subsections until \ref{47Inclusions} analyze the remaining cases in much the same fashion as in \cite{Pink2013b}. The last two subsections \ref{48Monodromy} and \ref{49MonodromyIncl} apply the preceding results to the actual geometric monodromy group $G^\geom$ of a quadratic morphism and deduce some consequences for the associated arithmetic monodromy group $G^\arith$.


\subsection{Groups associated to infinite postcritical orbits}
\label{41TheGroups}

Consider a set $X$ with a map $f\colon X\to X$. Consider a subset $C\subset X$ consisting of two distinct elements $p_0$ and $q_0$ which satisfy $f(p_0)\not=f(q_0)$. Assume that $P := \bigcup_{n\ge1}f^n(C) \subset X$ is infinite. In this subsection we study elements of $W$ which satisfy the same kind of recursion relations and product relation that hold for the generators of the geometric monodromy group of a quadratic morphism with an infinite postcritical orbit~$P$, but we do not assume that they come from an actual quadratic morphism.

Suppose that for every $p\in P$ we are given an element $b_p\in W$. Let $G$ denote the closure of the subgroup of $W$ that is generated by the $b_p$ for all $p\in P$.
We call the elements $b_p$ and the group $G$ \emph{weakly of type} $(X,f,C)$ if for every $p\in P$ we have
\UseTheoremCounterForNextEquation
\begin{equation}\label{4GenRecConj1}
b_p\ \sim\ \left\{\begin{array}{ll}
\ \ \ \sigma & \hbox{if $p=f(c)$ for some $c\in C\setminus P$,}\\[3pt]
(b_c,1)\,\sigma & \hbox{if $p=f(c)$ for some $c\in C\cap P$,}\\[3pt]
(b_q,1) & \hbox{if $p=f(q)$ for a unique $q\in P\setminus C$,} \\[3pt]
(b_q,b_{q'}) & \hbox{if $p=f(q)=f(q')$ for distinct $q,q'\in P\setminus C$.}
\end{array}\right\}
\end{equation}
We call the elements $b_p$ and the group $G$ \emph{strongly of type} $(X,f,C)$ if, in addition, the infinite product of all $b_p$ in some order converges to the identity element of~$W$. The main content of this article is a study of all groups that are strongly of type $(X,f,C)$.

\medskip
For this we first classify the combinatorial possibilities for $P$ together with the map $P\to P$ induced by~$f$. Abbreviate $p_n := f^n(p_0)$ and $q_n := f^n(q_0)$ for all $n\ge1$.

\begin{Prop}\label{41Class}
We have precisely one of the following cases:
\begin{enumerate}
\item[(a)] The elements $p_1,p_2,\ldots$ and $q_1,q_2,\ldots$ are all distinct.
\item[(b)] The elements $p_1,p_2,\ldots$ and $q_1,\ldots,q_r$ are all distinct, and $q_{r+1}=q_{s+1}$, for unique indices $r>s\ge0$.
\item[(b\kern1pt$'$\kern-2pt)] Same as (b) with the $p_n$ and $q_n$ interchanged.
\item[(c)] The elements $p_1,p_2,\ldots$ and $q_1,\ldots,q_s$ are all distinct, and $p_{r+1}=q_{s+1}$, for unique indices $r,s\ge0$ which are not both~$0$.
\end{enumerate}
\end{Prop}

\begin{Proof}
If the elements $q_1,q_2,\ldots$ are not all distinct, there exists a relation of the form $q_{r+1}=q_{s+1}$ for some $r>s\ge0$. Then $q_{r+i}=q_{s+i}$ for all $i\ge1$, and so $\{q_1,q_2,\ldots\} = \{q_1,\ldots,q_r\}$ is finite. The same remark applies with $p_n$ in place of~$q_n$. Since $P$ is infinite, this cannot occur for both the $p_n$ and the $q_n$. After possibly interchanging $p_0$ with $q_0$, which interchanges the cases (b) and (b$'$), we may thus without loss of generality assume that the elements $p_1,p_2,\ldots$ are all distinct.

Suppose in addition that the elements $q_1,q_2,\ldots$ are all distinct from $p_1,p_2,\ldots$. If they are also distinct from each other, we have the case (a). Otherwise there exists a unique smallest $r\ge1$ such that $q_{r+1}=q_{s+1}$ for some $s$ satisfying $r>s\ge0$. Then the elements $q_1,\ldots,q_r$ are distinct and $s$ is unique, so we have the case (b).

Suppose now that the elements $q_1,q_2,\ldots$ are not all distinct from $p_1,p_2,\ldots$. Then there exists a unique smallest $s\ge0$ such that $q_{s+1}=p_{r+1}$ for some $r\ge0$. This $r$ is then also unique. Moreover, any relation of the form $q_i=q_j$ for $1\le i<j\le s$ would imply that $p_{r+i}=q_{s+i}=q_{s+j}=p_{r+j}$, contradicting the assumption that the elements $p_1,p_2,\ldots$ are distinct. Thus the elements $q_1,\ldots,q_s$ are all distinct and distinct from $p_1,p_2,\ldots$. Finally, since by assumption $p_1\not=q_1$, we cannot have $r=s=0$. Thus we have the case (c).
\end{Proof}

\medskip
Here is a graphic depiction in the respective cases of the set $P$, where the map $P\to P$ induced by $f$ is represented by arrows, and the `entry points' $p_1$ and~$q_1$ are specially marked:
$$\vtop{\hsize=200pt
\hbox{Case (a)}
\vskip5pt\noindent
\fbox{\ $
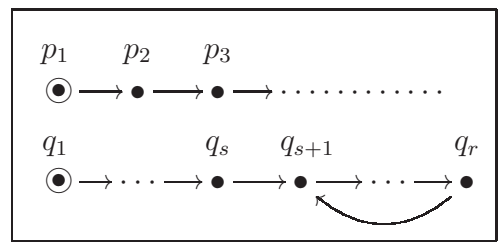
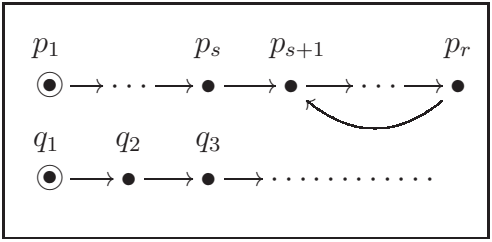
\xymatrix@R-26pt@C-12pt{
p_1 & p_2 & p_3 & \phantom{p_n} \raise3pt\hbox{$\mathstrut$}\\
\circbullet\ar[r]&\bullet\ar[r]&\bullet\ar[r]&\cdots \\
\raise-2pt\hbox{$q_1$} & \raise-2pt\hbox{$q_2$}  & \raise-2pt\hbox{$q_3$} & \phantom{q_n}\raise3pt\hbox{$\mathstrut$} \\
\smash{\circbullet}\ar[r]&\bullet\ar[r]&\bullet\ar[r]&\cdots
{\raise-7pt\hbox{\vphantom{.}}} \\
}$}
}\qquad
\vtop{\hsize=180pt
\hbox{Case (b) with $r>s\ge0$}
\vskip5pt\noindent
\fbox{\ $\xymatrix@R-26pt@C-12pt{
p_1 & p_2 & p_3 & \phantom{p_n} \raise3pt\hbox{$\mathstrut$} &\\
\circbullet\ar[r]&\bullet\ar[r]&\bullet\ar[r]&\cdots\rlap{$\cdots\cdots\cdots$}  \\
\raise-2pt\hbox{$q_1$} & \phantom{\raise-2pt\hbox{$q_n$}} & \raise-2pt\hbox{$q_s$} & 
\kern-5pt \raise-2pt\hbox{$q_{s+1}$} \kern-10pt & \phantom{\raise-2pt\hbox{$q_n$}} & 
\raise-2pt\hbox{$q_r$} \raise3pt\hbox{$\mathstrut$}  \\ 
\circbullet \ar[r]&\cdots\ar[r]&\bullet\ar[r]&\bullet\ar[r]&\cdots\ar[r]&\bullet\ar@/^16pt/[ll] \\
&&&&& {\raise-3pt\hbox{\vphantom{.}}} \\
}$}
}$$
$$\vtop{\hsize=200pt
\hbox{Case (b$\kern1pt'$) with $r>s\ge0$}
\vskip5pt\noindent
\fbox{\ $\xymatrix@R-26pt@C-12pt{
\raise0pt\hbox{$p_1$} & \phantom{\raise0pt\hbox{$p_n$}} & \raise0pt\hbox{$p_s$} & 
\kern-5pt \raise0pt\hbox{$p_{s+1}$} \kern-10pt & \phantom{\raise0pt\hbox{$p_n$}} & 
\raise0pt\hbox{$p_r$} \raise3pt\hbox{$\mathstrut$}  \\ 
\circbullet \ar[r]&\cdots\ar[r]&\bullet\ar[r]&\bullet\ar[r]&\cdots\ar[r]&\bullet\ar@/^16pt/[ll] \\
\raise-2pt\hbox{$q_1$} & \raise-2pt\hbox{$q_2$} & \raise-2pt\hbox{$q_3$} & \phantom{q_n} \raise12pt\hbox{\vphantom{.}} &\\
\circbullet\ar[r]&\bullet\ar[r]&\bullet\ar[r]&\cdots\rlap{$\cdots\cdots\cdots$}  \\
&&&&& {\raise-3pt\hbox{\vphantom{.}}} \\
}$} 
}\qquad
\vtop{\hsize=180pt
\hbox{Case (c) with $r,s\ge0$, not both $0$}
\vskip5pt\noindent
\fbox{\ $
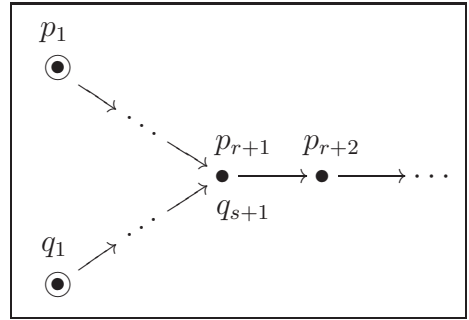
\xymatrix@R-21pt@C-10pt{
\vphantom{x}\myothersmash{-1pt}{p_1} &&&& \\
\circbullet \ar[dr] &&&& \\
&\raisebox{-1.5pt}{$\ddots$}\ar[dr]& 
\kern-3pt\myothersmash{-6pt}{\ p_{r+1}}\kern-15pt
& \kern7pt\myothersmash{-6pt}{p_{r+2}}\kern0pt & \\
&&\bullet\ar[r]&\bullet\ar[r]&\cdots\\
\myothersmash{-6pt}{q_1} &\raisebox{-1.5pt}{\reflectbox{$\ddots$}}\ar[ur]&
\kern-3pt\myothersmash{8pt}{\ q_{s+1}}\kern-15pt &&\\
\circbullet \ar[ur]&&&&  {\raise-3pt\hbox{\vphantom{.}}}\  \\
}$} 
}$$

\medskip
In most cases any group  of type $(X,f,C)$ is equal to~$W$. Namely:

\begin{Prop}\label{41LargeAB}
In the cases (a) and (b) and (b\kern1pt$'$\kern-2pt) of Proposition \ref{41Class}, any group that is weakly of type $(X,f,C)$ is equal to~$W$.
\end{Prop}

\begin{Proof}
By symmetry the case (b$'$) reduces to the case (b). In the cases (a) and (b) the elements $p_1,p_2,\ldots$ are all distinct and not of the form $f(q_i)$ for any $i\ge0$. Thus the relations (\ref{4GenRecConj1}) imply that 
$$\biggl\{\begin{array}{ccll}
b_{p_1}     \!\!\!&\sim&\!\!\! \sigma & \hbox{\ and}\\[3pt]
b_{p_i}     \!\!\!&\sim&\!\!\! (b_{p_{i-1}},1) & \hbox{\ for all $i>1$.} \\[3pt]
\end{array}\biggr\}$$
As signs are invariant under conjugation, it follows that $\sgn_1(b_{p_1})=-1$ and $\sgn_n(b_{p_1}) =\sgn_1(b_{p_i})=\nolinebreak 1$ and $\sgn_n(b_{p_i}) =\sgn_{n-1}(b_{p_{i-1}})$ for all $n,i>1$. By induction this implies that 
$$\sgn_n(b_{p_i}) \ =\ 
(-1)^{\delta_{n,i}}\quad\hbox{(Kronecker delta)}$$
for all $n,i\ge1$. With \cite[Prop.\ \ref{15SignGnSurj} (b)]{Pink2013b} we conclude that $G=W$. 
\end{Proof}

\begin{Prop}\label{41LargeC}
In the case (c) of Proposition \ref{41Class} with $r\not=s$, any group that is weakly of type $(X,f,C)$ is equal to~$W$.
\end{Prop}

\begin{Proof}
By symmetry we may without loss of generality assume that $r>s$. Then the relations (\ref{4GenRecConj1}) imply that
$$\left\{\begin{array}{lcll}
b_{p_1}     \!\!\!&\sim&\!\!\! \sigma & \\[3pt]
b_{p_i}     \!\!\!&\sim&\!\!\! (b_{p_{i-1}},1) & \hbox{\ for all $i>1$ with $i\not=r+1$,}
 \\[3pt]
b_{p_{r+1}} \!\!\!&\sim&\!\!\! (b_{p_r},1)\,\sigma & \hbox{\ if $s=0$,} \\[3pt]
b_{p_{r+1}} \!\!\!&\sim&\!\!\! (b_{p_r},b_{q_s}) & \hbox{\ if $s>0$,} \\[3pt]
b_{q_1}     \!\!\!&\sim&\!\!\! \sigma & \hbox{\ if $s>0$,} \\[3pt]
b_{q_j}     \!\!\!&\sim&\!\!\! (b_{q_{j-1}},1) & \hbox{\ for all $1<j\le s$.} \\[3pt]
\end{array}\right\}$$
As in the proof of Proposition \ref{41LargeAB} we deduce from this by induction that for all $n\ge1$
$$\begin{array}{ll}
\sgn_n(b_{p_i}) = (-1)^{\delta_{n,i}} & \hbox{for all $1\le i\le r$, and}\\[5pt]
\sgn_n(b_{q_j}) = (-1)^{\delta_{n,j}} & \hbox{for all $1\le j\le s$.}
\end{array}$$
Moreover we claim that for all $n\ge1$ and all $i\ge r+1$ we have
$$\sgn_n(b_{p_i}) = (-1)^{\delta_{n,i}+\delta_{n,i-r+s}}.\ \qquad\qquad\qquad$$
Indeed, for $i=r+1$ one sees this by going doggedly through all the possible cases, which we leave to the persistent reader. For $i>r+1$ it follows easily by induction using the fact that $b_{p_i} \sim (b_{p_{i-1}},1)$.

Now define elements $c_i\in G$ for all $i\ge1$ by setting recursively
$$c_i \ :=\ 
\biggl\{\!\begin{array}{ll}
b_{p_i} & \hbox{if $i\le r$,}\\[3pt]
b_{p_i}\kern1pt c_{i-r+s} & \hbox{if $i>r$.}
\end{array}$$
Using the above formulas for all $\sgn_n(b_{p_i})$, by induction on $i$ one easily shows that
$$\sgn_n(c_i) = (-1)^{\delta_{n,i}}$$
for all $n,i\ge1$. With \cite[Prop.\ \ref{15SignGnSurj} (b)]{Pink2013b} we conclude that $G=W$. 
\end{Proof}

\medskip
The only remaining case is the case (c) of Proposition \ref{41Class} with $r=s\ge1$. In this case the relations (\ref{4GenRecConj1}) are equivalent to
\UseTheoremCounterForNextEquation
\begin{equation}\label{4GenRecConj2}
\left\{\begin{array}{lcll}
b_{p_1}     \!\!\!&\sim&\!\!\! \sigma & \\[3pt]
b_{p_i}     \!\!\!&\sim&\!\!\! (b_{p_{i-1}},1) & \hbox{\ for all $i>1$ with $i\not=r+1$,}
 \\[3pt]
b_{p_{r+1}} \!\!\!&\sim&\!\!\! (b_{p_r},b_{q_r}) & \\[3pt]
b_{q_1}     \!\!\!&\sim&\!\!\! \sigma & \\[3pt]
b_{q_j}     \!\!\!&\sim&\!\!\! (b_{q_{j-1}},1) & \hbox{\ for all $1<j\le r$.} \\[3pt]
\end{array}\right\}
\end{equation}
We will not say anything about groups which are only weakly of this type, because without the product relation for the generators they do not seem rigid enough. In the following subsections we will first study a single group that is strongly of this type, then we will show that any other is conjugate to it, and then we will study the chosen group some more. By combining Propositions \ref{41LargeAB} and \ref{41LargeC} with Theorem \ref{44SemiRigid} below we obtain:

\begin{Thm}\label{41All}
All subgroups which are strongly of type $(X,f,C)$ are conjugate under~$W$.
\end{Thm}


\subsection{Setup and basic properties}
\label{42Basic}

Consider an integer $r\ge1$. Using \cite[Prop. \ref{15RecRelsProp}]{Pink2013b} we define elements $a_1,a_2,\ldots\in W$ and $b_1,\ldots,b_r\in W$ by the recursion relations
\UseTheoremCounterForNextEquation
\begin{equation}\label{4RecRels}
\left\{\begin{array}{lcll}
a_1     \!\!\!&=&\!\!\! \sigma & \\[3pt]
a_{r+1} \!\!\!&=&\!\!\! (a_r,b_r) & \\[3pt]
a_i     \!\!\!&=&\!\!\! (a_{i-1},1) & \hbox{\ for all $i>1$ with $i\not=r+1$,}
 \\[3pt]
b_1     \!\!\!&=&\!\!\! (b_r,b_r^{-1})\,\sigma & \\[3pt]
b_j     \!\!\!&=&\!\!\! (b_rb_{j-1}b_r^{-1},1) & \hbox{\ for all $1<j\le r$.} \\[3pt]
\end{array}\right\}
\end{equation}
Let $G(r)$ denote the closure of the subgroup generated by all these elements. Since $r$ will be fixed until the end of Subsection \ref{45Normalizer}, we abbreviate $G := G(r)$ until then.

\begin{Prop}\label{42GenOrderTriv}
\begin{enumerate}
\item[(a)] Every generator $a_i$ and $b_j$ has order~$2$.
\item[(b)] 
For all $i>n\ge0$ we have $a_i|_{T_n} = 1$, and for all $r\ge j>n\ge0$ we have $b_j|_{T_n} = 1$.
\end{enumerate}
\end{Prop}

\begin{Proof}
For $a_1$ and $b_1$ both statements follow directly from the recursion relations (\ref{4RecRels}), for the other generators they follow by induction on $i$ and~$j$.
\end{Proof}

\begin{Prop}\label{42Product}
We have \ $b_1\cdots b_r = a_1a_2a_3\cdots$.
\end{Prop}

\begin{Proof}
Proposition \ref{42GenOrderTriv} (b) implies that the infinite product on the right hand side converges in~$W$. Using the recursion relations (\ref{4RecRels}) and the fact that $b_r^{-1}=b_r$ we calculate
\begin{eqnarray*} 
\alpha\; :=\; a_1a_2a_3\cdots 
&\!\!=\!\!& \sigma\, (a_1,1)\cdots(a_{r-1},1)\, (a_r,b_r)
\,(a_{r+1},1)\, (a_{r+2},1)\cdots \\
&\!\!=\!\!& \sigma\, (a_1\cdots a_{r-1}a_ra_{r+1}a_{r+2}\cdots, b_r) \\
&\!\!=\!\!& \sigma\, (\alpha, b_r) \hskip180pt \hbox{and} \\
\beta\ :=\; b_1\cdots b_r\;\ \ 
&\!\!=\!\!& (b_r,b_r^{-1})\,\sigma\, (b_rb_1b_r^{-1},1) \cdots (b_rb_{r-1}b_r^{-1},1) \\
&\!\!=\!\!& \sigma\,(b_r^{-1},b_r)\, (b_rb_1\cdots b_{r-1}b_r^{-1},1) \\
&\!\!=\!\!& \sigma\,(b_1\cdots b_{r-1}b_r^{-1},b_r) \\
&\!\!=\!\!& \sigma\, (\beta, b_r).
\end{eqnarray*} 
Thus
$$\gamma\ :=\ \alpha^{-1}\beta
\ =\ (\alpha^{-1}, b_r^{-1})\,\sigma^{-1}\,\sigma\,(\beta,b_r) 
\ =\ (\alpha^{-1}\beta,1)
\ =\ (\gamma,1).$$
By \cite[Prop.\ \ref{15RecTriv}]{Pink2013b} it follows that $\gamma=1$ and hence $\beta=\alpha$, as desired.
\end{Proof}

\begin{Prop}\label{42StrongType}
The generators $a_i$ and $b_j$ and the group $G$ are strongly of the type (c) of Proposition \ref{41Class} with $s=r$. 
\end{Prop}

\begin{Proof}
Comparison of (\ref{4RecRels}) with (\ref{4GenRecConj2}) shows that the generators $a_i$ and~$b_j$, though renamed, are weakly of this type.
On the other hand, since all $b_j$ have order~$2$, the equation in Proposition \ref{42Product} is equivalent to $b_r\cdots b_1a_1a_2a_3\cdots = 1$.
Thus the generators and hence the group $G$ are strongly of the required type.
\end{Proof}

\medskip
It turns out that the recursion relations (\ref{4RecRels}) are only one choice among many yielding a group that is strongly of the above type. In Subsection \ref{44Conjugacy} we will show any other choice yields a group which is conjugate to~$G$.

\begin{Prop}\label{41GenSigns}
For all $n,i\ge1$ we have
$$\left\{\!\begin{array}{ll}
\sgn_n(a_i) = \sgn_n(b_i) = -1 & \hbox{if $n=i\le r$,}\\[3pt]
\sgn_n(a_i) = \sgn_n(b_i) = 1 & \hbox{if $n\not=i\le r$,}\\[3pt]
\sgn_n(a_i) = 1 & \hbox{if $i>r$.}
\end{array}\right\}$$
\end{Prop}

\begin{Proof}
The recursion relations (\ref{4RecRels}) imply this directly if $n=1$ or $i=1$. If $n>1$ and $1<i\le r$, they show that $\sgn_n(a_i) = \sgn_{n-1}(a_{i-1})$ and $\sgn_n(b_i) = \sgn_{n-1}(b_{i-1})$, so the desired equations for all $i\le r$ follow by induction. In particular, this shows that $\sgn_n(a_r)=\sgn_n(b_r)$ for all $n\ge1$. Thus the recursion relation implies that $\sgn_n(a_{r+1}) = \sgn_n((a_r,b_r)) = 1$ for all $n\ge1$. By induction this in turn implies that $\sgn_n(a_i) = 1$ whenever $i\ge r+1$, and we are done.
\end{Proof}

\begin{Prop}\label{41GnAll}
For any $n\ge0$ we have $G_n=W_n$ if and only if $n\le r$.
\end{Prop}

\begin{Proof}
Direct consequence of \cite[Prop.\ \ref{15SignGnSurj} (a)]{Pink2013b} and Proposition \ref{41GenSigns}.
\end{Proof}


\subsection{Useful subgroups}
\label{42Subgroups}

First we consider the subgroup of index $2$ which acts trivially on level~$1$:
\UseTheoremCounterForNextEquation
\begin{equation}\label{23G1Def}
G^1\ :=\ G\cap(W\times W).
\end{equation}
Let $\proj_1$ and $\proj_2\colon W\times W\to W$ denote the two projections. We have the following self-similarity properties:

\begin{Prop}\label{42G1Prop}
\begin{enumerate}
\item[(a)] $G\subset (G\times G)\rtimes \langle\sigma\rangle$.
\item[(b)] $G^1\subset G\times G$.
\item[(c)] $\proj_1(G^1)=\proj_2(G^1)=G$.
\end{enumerate}
\end{Prop}

\begin{Proof}
By the recursion relations (\ref{4RecRels}), the generators of $G$ lie in $(G\times G)\rtimes \langle\sigma\rangle$, hence so does~$G$, proving (a). Also, (a) directly implies (b). 

By (\ref{4RecRels}) the subgroup $G^1$ contains the elements $(a_{i-1},1)$ for all $i>1$ with $i\not=r+1$, or equivalently the elements $(a_i,1)$ for all $i\ge1$ with $i\not=r$. Since $G^1$ also contains $(a_r,b_r)$, it follows that $\proj_1(G^1)$ contains $a_i$ for all $i\ge1$. Similarly, since $G^1$ contains the elements $(b_rb_{j-1}b_r^{-1},1)$ for all $1<j\le r$, it follows that $\proj_1(G^1)$ contains $b_rb_jb_r^{-1}$ for all $1\le j<r$. But $G^1$ also contains the element $b_1a_1 = (b_r,b_r^{-1})$; hence 
$\proj_1(G^1)$ contains~$b_r$. Thus it contains $b_j$ for all $1\le j\le r$, and is therefore equal to~$G$.

Finally, conjugation by the element $a_1=\sigma \in G$ interchanges the two factors of $W\times W$ and normalizes~$G^1$; hence also $\proj_2(G^1)=G$, proving (c).
\end{Proof}

\begin{Prop}\label{42LevelTrans}
The group $G$ acts transitively on the level $n$ of $T$ for every $n\ge1$.
\end{Prop}

\begin{Proof}
As $G$ contains $a_1=\sigma$, it acts transitively on level~$1$. If it acts transitively on level~$n$, Proposition \ref{42G1Prop} (c) implies that $G^1$ acts transitively on the subset of level $n+1$ lying over any fixed vertex of level~$1$. As $G$ permutes the level $1$ transitively, it follows that $G$ acts transitively on level $n+1$. Thus the proposition follows by induction.
\end{Proof}

\medskip
Next we consider the following normal subgroup of~$G$:
\UseTheoremCounterForNextEquation
\begin{equation}\label{42HDef}
H\ :=\ \biggl\{\begin{array}{l}
\hbox{closure of the subgroup generated by all $G$-conjugates} \\[3pt]
\hbox{of $a_i$ for all $i\not=r$,\ \ of $b_j$ for all $j<r$,\ \ and of $a_rb_r$.} \\[3pt]
\end{array}\biggr\}_.
\end{equation}

\begin{Prop}\label{42HSemiDirect}
We have:
\begin{enumerate}
\item[(a)] $H=\Ker(\sgn_r|_G)$, 
\item[(b)] $G = H \rtimes \langle b_r\rangle$,
\item[(c)] $G = (H\!\times\! H) \rtimes 
\bigl( \langle(b_r,b_r)\rangle \!\times\! \langle\sigma\rangle \bigr)$,
\item[(d)] $G^1 = (H\!\times\! H) \rtimes \langle (b_r,b_r) \rangle$,
\end{enumerate}
where $\langle b_r\rangle$ and $\langle(b_r,b_r)\rangle$ and $\langle\sigma\rangle$ are cyclic groups of order~$2$.
\end{Prop}

\begin{Proof}
As $b_r$ is an element of order $2$ with $\sgn_r(b_r)=-1$, we have $G = {\Ker(\sgn_r|_G) \rtimes \langle b_r\rangle}$.
The results on signs in Proposition \ref{41GenSigns} also show that all generators of~$H$, and hence $H$ as well, are contained in $\Ker(\sgn_r|_G)$.
On the other hand the definitions of $G$ and $H$ imply that the factor group $G/H$ is topologically generated by the image of~$b_r$. Therefore $[G:H]\le 2$, which together leaves only the possibility $H=\Ker(\sgn_r|_G)$. This proves (a) and at the same time (b).

Next consider the subgroup $H':=\{x\in G\mid(x,1)\in G\}$. For any $u\in G$, by Proposition \ref{42G1Prop} (c) there exists $v\in G$ such that $(u,v)\in G$. For any $x\in H'$ we then have $(uxu^{-1},1) = (u,v)\,(x,1)\,(u,v)^{-1} \in G$ and hence $uxu^{-1}\in H'$. Thus $H'$ is a normal subgroup of~$G$.

By (\ref{4RecRels}) the subgroup $G^1$ contains the elements $(a_{i-1},1)$ for all $i>1$ with $i\not=r+1$, that is to say, the elements $(a_i,1)$ for all $i\ge1$ with $i\not=r$. Similarly $G^1$ contains the elements $(b_rb_jb_r^{-1},1)$ for all $1\le j<r$. It also contains the element $a_{r+1}b_1a_1 = (a_r,b_r)\,(b_r,b_r^{-1})\,\sigma\,\sigma = (a_rb_r,1)$. Thus $H'$ contains the elements $a_i$ for all $i\not=r$ and $b_rb_jb_r^{-1}$ for all $j<r$ and $a_rb_r$.
As $H'$ is normal in~$G$, it thus contains all generators of $H$ and therefore $H$ itself. This shows that $H\times 1\subset G$.

Since $H$ is normal in~$G$, Proposition \ref{42G1Prop} (c) implies in the same way as above that $H\times 1$ is normal in~$G^1$. As its conjugate by $a_1=\sigma$ is $1\times H$, it follows that $H\times H$ is a normal subgroup of~$G$.

Moreover (\ref{4RecRels}) shows that $H\times H$ contains all generators $a_i,b_j$ of $G$ except $a_1,b_1,a_{r+1}$. Thus the factor group $G/(H\times H)$ is topologically generated by the images of $a_1,b_1,a_{r+1}$. By the product relation from Proposition \ref{42Product} we can drop $a_{r+1}$ and deduce that $G/(H\times H)$ is topologically generated by the images of $a_1=\sigma$ and $b_1a_1 = (b_r,b_r^{-1})\,\sigma\,\sigma = (b_r,b_r)$ alone. As these two elements commute and have order~$2$, it follows that $G$ is the almost semidirect product 
$G = (H\!\times\! H) \cdot \bigl( \langle(b_r,b_r)\rangle \!\times\! \langle\sigma\rangle \bigr)$. Intersecting with $W\times W$ we deduce that 
$G^1$ is the almost semidirect product $G^1 = (H\!\times\! H) \cdot \langle(b_r,b_r)\rangle$. But since $b_r\not\in H$ by (b), this almost semidirect product is a true semidirect product, proving (d). Since $\sigma\not\in G^1$, this in turn implies (c), and the proposition is proved.
\end{Proof}

\medskip
We will also need an analogue of Proposition \ref{42HSemiDirect} on finite levels. For any subgroup $X\subset W$ we let $X_n$ denote its image in~$W_n$. 

\begin{Prop}\label{42nHSemiDirect}
For all $n\ge r$ we have
\begin{enumerate}
\item[(a)] $b_r|_{T_n}\not\in H_n$,
\item[(b)] $G_n = H_n \rtimes \langle b_r|_{T_n}\rangle$,
\item[(c)] $G_{n+1} = (H_n\!\times\! H_n) \rtimes 
\bigl( \langle(b_r|_{T_n},b_r|_{T_n})\rangle \!\times\! \langle\sigma\rangle \bigr)$, and 
\item[(d)] $G^1_{n+1} = (H_n\!\times\! H_n) \rtimes \langle (b_r|_{T_n},b_r|_{T_n}) \rangle$,
\end{enumerate}
where $\langle b_r|_{T_n}\rangle$ and $\langle(b_r|_{T_n},b_r|_{T_n})\rangle$ and $\langle\sigma\rangle \subset W_{n+1}$ are cyclic groups of order~$2$.
\end{Prop}

\begin{Proof}
For all $n\ge r$ the homomorphism $\sgn_r$ factors through a homomorphism $W_n\to\{\pm1\}$. Proposition \ref{42HSemiDirect} (a) implies that the latter is trivial on~$H_n$, but nontrivial on~$b_r|_{T_n}$. This implies (a). The remaining assertions follow from (a) and the corresponding assertions in \ref{42HSemiDirect}.
\end{Proof}



\subsection{Size}
\label{43Size}

\begin{Prop}\label{43GnOrder}
For all $n\ge0$ we have 
$$\log_2|G_n|\ =
\biggl\{\begin{array}{ll}
2^n-1         & \hbox{if $n\le r$,} \\[3pt]
2^n - 2^{n-r} & \hbox{if $n\ge r$.}
\end{array}$$
\end{Prop}

\begin{Proof}
For $n\le r$ this results from Proposition \ref{41GnAll} and the formula (\ref{12WnOrder}) for $\log_2|W_n|$ from \cite{Pink2013b}. For $n\ge r$ we use Proposition \ref{42nHSemiDirect} (c) and (b) to calculate $|G_{n+1}| = 4\cdot|H_n|^2 = |G_n|^2$.
Thus $\log_2|G_{n+1}| = 2\cdot\log_2|G_n|$, which by induction implies that
$\log_2|G_n| = {2^{n-r}\cdot\log_2|G_r|} = 
2^{n-r}\cdot(2^r-1) = 2^n-2^{n-r}$ for all $n\ge r$, as desired.
\end{Proof}

\begin{Thm}\label{43Hausdorff}
The Hausdorff dimension of $G$ exists and is $1-2^{-r}$.
\end{Thm}

\begin{Proof}
The Hausdorff dimension of $G$ is defined as the limit of $\frac{\log_2|G_n|}{2^n-1}$ for $n\to\infty$. From Proposition \ref{43GnOrder} we find the desired value.
\end{Proof}


\subsection{Conjugacy of generators}
\label{44Conjugacy}

Recall that $\sim$ means conjugacy under~$W$. 

\begin{Thm}\label{44SemiRigid}
Consider elements $a'_i$, $b'_j\in W$ for all $i\ge1$ and $1\le j\le r$, which satisfy:
\begin{enumerate}
\item[(a)] 
$\left\{\begin{array}{lcll}
a'_1     \!\!\!&\sim&\!\!\! \sigma & \\[3pt]
a'_{r+1} \!\!\!&\sim&\!\!\! (a'_r,b'_r) & \\[3pt]
a'_i     \!\!\!&\sim&\!\!\! (a'_{i-1},1) & \hbox{\ for all $i>1$ with $i\not=r+1$,}
 \\[3pt]
b'_1     \!\!\!&\sim&\!\!\! \sigma & \\[3pt]
b'_j     \!\!\!&\sim&\!\!\! (b'_{j-1},1) & \hbox{\ for all $1<j\le r$,} \\[3pt]
\end{array}\right\}$\ \ and
\item[(b)] their infinite product in some order converges to the identity element of~$W$. 
\end{enumerate}
Then there exists $w\in W$ such that for all $i$ and~$j$, the element $w^{-1}a'_iw$ is conjugate to $a_i$ under~$G$, and $w^{-1}b'_jw$ is conjugate to $b_j$ under~$G$. 
Moreover, for any such $w$ the closure of the subgroup of $W$ generated by all $a'_i$ and $b'_j$ is $wGw^{-1}$.
\end{Thm}

Comparison with (\ref{4GenRecConj2}) shows that the assumptions mean that the generators $a'_i$ and~$b'_j$, though renamed, and the closure of the subgroup generated by them, are strongly of the type (c) of Proposition \ref{41Class} with $s=r$. In particular Theorem \ref{44SemiRigid} implies that any two subgroups which are strongly of this type are conjugate under~$W$.

\bigskip
\begin{Proof}
The second statement follows from the first by \cite[Lemma \ref{13ConjGen}]{Pink2013b}. Also, by \cite[Prop.\ \ref{15RecConj=Conj}]{Pink2013b}, the assumption (a) implies that $a'_i\sim a_i$ and $b'_j\sim b_j$ for all $i$ and~$j$. By Proposition \ref{42GenOrderTriv} (b) we therefore have $a'_i|_{T_n} = 1$ for all $i>n$. In particular, the product of the restrictions $a'_i|_{T_n}$ and $b'_j|_{T_n}$ in any order is essentially a finite product. Thus the first statement follows by taking the limit over $n$ of the following assertion for all $n\ge0$: 
\begin{enumerate}
\item[($*_n$)] For any elements $a'_i$, $b'_j\in W$ satisfying $a'_i\sim a_i$ and $b'_j\sim b_j$, such that the product of all $a'_i|_{T_n}$ and $b'_j|_{T_n}$ in some order is equal to $1$, there exists $w\in W$ such that for all $i$ and~$j$, the element $w^{-1}a'_iw|_{T_n}$ is conjugate to $a_i|_{T_n}$ under~$G_n$, and $w^{-1}b'_jw|_{T_n}$ is conjugate to $b_j|_{T_n}$ under~$G_n$. 
\end{enumerate}
This is trivial for $n=0$, so assume that $n>0$ and that ($*_{n-1}$) is true. Take elements $a'_i$, $b'_j\in W$ satisfying the assumptions in ($*_n$). Then in particular $a'_1\sim a_1=\sigma$. After conjugating everything by the same element of $W$ we may therefore without loss of generality assume that $a_1=\sigma$. Since $b'_1\sim b_1\sim\sigma$, we also have $b'_1=(c,c^{-1})\,\sigma$ for some element $c\in W$, with which we will deal later. For the remaining elements we have
$$\left\{\begin{array}{lll}
a'_{r+1} &\!\!\!\sim\  a_{r+1} \,=\, (a_r,b_r), & \\[3pt]
a'_i &\!\!\!\sim\ a_i \,=\, (a_{i-1},1) & \hbox{for all $i\not=1,r+1$,} \\[3pt]
b'_j &\!\!\!\sim\ b_j \,=\, (b_rb_{j-1}b_r^{-1},1) \,\sim\, (b_{j-1},1) & \hbox{for all $1<j\le r$.}  
\end{array}\right\}$$
By basic properties of conjugacy in $W$ (see \cite[Lemma \ref{13ConjEquiv}]{Pink2013b}) this means that
\UseTheoremCounterForNextEquation
\begin{equation}\label{44SemiRigidLem1}
\left\{\begin{array}{lcl}
a'_{r+1} &\!\!\!=\!\!& \hbox{\,$(a''_r,b''_r)$ \:or\, $(b''_r,a''_r)$}, \\[3pt]
a'_i &\!\!\!=\!\!& \hbox{$(a''_{i-1},1)$ or $(1,a''_{i-1})$} \ \ \hbox{for all $i\not=1,r+1$,} \\[3pt]
b'_j &\!\!\!=\!\!& \hbox{$(b''_{j-1},1)$ or $(1,b''_{j-1})$} \ \ \hbox{for all $1<j\le r$,}  
\end{array}\right\}
\end{equation}
for elements $a''_i$, $b''_j\in W$ satisfying $a''_i\sim a_i$ and $b''_j\sim b_j$ for all $i\ge1$ and $1\le j\le r$.

\medskip
Next, by assumption and (\ref{44SemiRigidLem1}) the product of the elements 
$$\left\{\begin{array}{l}
\,\sigma|_{T_n}, \\[3pt]
(c,c^{-1})\,\sigma|_{T_n}, \\[3pt]
\hbox{$(a''_r,b''_r)|_{T_n}$ \,or\, $(b''_r,a''_r)|_{T_n}$}, \\[3pt]
\hbox{$(a''_{i-1},1)|_{T_n}$ or $(1,a''_{i-1})|_{T_n}$} \ \ \hbox{for all $i\not=1,r+1$,} \\[3pt]
\hbox{$(b''_{j-1},1)|_{T_n}$ or $(1,b''_{j-1})|_{T_n}$} \ \ \hbox{for all $1<j\le r$,}  
\end{array}\right\}$$
in some order is equal to~$1$. Solving this equation for the element $(c,c^{-1})\,\sigma|_{T_n}$, which is equal to its own inverse, shows that this element is equal to the product of the others in some order. In the resulting equation we can move the factor $\sigma|_{T_n}$ to the right by interchanging the entries of the intervening factors. After multiplying by $\sigma^{-1}|_{T_n}$ from the right hand side we can deduce that $(c,c^{-1})|_{T_n}$ is equal to the product of the elements
$$\left\{\begin{array}{l}
\hbox{$(a''_r,b''_r)|_{T_n}$ \,or\, $(b''_r,a''_r)|_{T_n}$}, \\[3pt]
\hbox{$(a''_{i-1},1)|_{T_n}$ or $(1,a''_{i-1})|_{T_n}$} \ \ \hbox{for all $i\not=1,r+1$,} \\[3pt]
\hbox{$(b''_{j-1},1)|_{T_n}$ or $(1,b''_{j-1})|_{T_n}$} \ \ \hbox{for all $1<j\le r$,}  
\end{array}\right\}$$
in some order (in general with other cases than before). In other words $(c|_{T_{n-1}},c^{-1}|_{T_{n-1}})$ is the product of the elements
\UseTheoremCounterForNextEquation
\begin{equation}\label{44SemiRigidLem2}
\left\{\begin{array}{l}
\hbox{$(a''_r|_{T_{n-1}},b''_r|_{T_{n-1}})$ 
\,or\, $(b''_r|_{T_{n-1}},a''_r|_{T_{n-1}})$}, \\[3pt]
\hbox{$(a''_i|_{T_{n-1}},1)$ or $(1,a''_i|_{T_{n-1}})$} \ \ \hbox{for all $i\not=r$,} \\[3pt]
\hbox{$(b''_j|_{T_{n-1}},1)$ or $(1,b''_j|_{T_{n-1}})$} \ \ \hbox{for all $1\le j<r$,}  
\end{array}\right\}
\end{equation}
in some order. Therefore $1 = c|_{T_{n-1}} \cdot c^{-1}|_{T_{n-1}}$ is the product of the elements
$$\left\{\begin{array}{l}
a''_r|_{T_{n-1}},\ b''_r|_{T_{n-1}}, \\[3pt]
a''_i|_{T_{n-1}} \ \ \hbox{for all $i\not=r$,} \\[3pt]
b''_j|_{T_{n-1}} \ \ \hbox{for all $1\le j<r$,}  
\end{array}\right\}$$
in some order. But this means that the product of all $a''_i|_{T_{n-1}}$ and $b''_j|_{T_{n-1}}$ in some order is equal to~$1$.

\medskip
We can thus apply the induction hypothesis ($*_{n-1}$) to the elements $a''_i$ and~$b''_j$, finding an element $u\in W$ such that for all $i$ and~$j$, the element $u^{-1}a''_iu|_{T_{n-1}}$ is conjugate to $a_i|_{T_{n-1}}$ under~$G_{n-1}$, and $u^{-1}b''_ju|_{T_{n-1}}$ is conjugate to $b_j|_{T_{n-1}}$ under~$G_{n-1}$. We then claim that $w := (u,u)\in W$ has the desired property in ($*_n$), that is to say, for all $i$ and~$j$ the element $w^{-1}a'_iw|_{T_n}$ is conjugate to $a_i|_{T_n}$ under~$G_n$, and $w^{-1}b'_jw|_{T_n}$ is conjugate to $b_j|_{T_n}$ under~$G_n$. 

\medskip
To see this for $w^{-1}a'_{r+1}w|_{T_n}$, choose elements $x,y\in G_{n-1}$ such that
$$\left\{\begin{array}{l}
u^{-1}a''_ru|_{T_{n-1}} \;=\; x\,(a_r|_{T_{n-1}})\,x^{-1}\ \ \hbox{and} \\[3pt]
u^{-1}b''_ru|_{T_{n-1}} \;=\; y\,(b_r|_{T_{n-1}})\,y^{-1}.
\end{array}\right\}$$
After possibly replacing $y$ by $y\,(b_r|_{T_{n-1}})$, which does not change the second equation, by Proposition \ref{42HSemiDirect} (b) we may without loss of generality assume that $H_{n-1} x=H_{n-1} y$. By Proposition \ref{42HSemiDirect} (d) this implies that $(x,y)$ lies in $G^1_n$. On the other hand (\ref{44SemiRigidLem1}) says that $a'_{r+1} = \sigma^\lambda\,(a''_r,b''_r)\,\sigma^{-\lambda}$ for some $\lambda\in\{0,1\}$. With $z := \sigma^\lambda\,(x,y) \in G_n$ we deduce that 
\begin{eqnarray*}
w^{-1}a'_{r+1}w|_{T_n}
&\!\!\!=\!\!\!& (u,u)^{-1}\,\sigma^\lambda (a''_r,b''_r)\,\sigma^{-\lambda}\, (u,u)|_{T_n} \\
&\!\!\!=\!\!\!& \sigma^\lambda\, (u^{-1}a''_ru|_{T_{n-1}}, u^{-1}b''_ru|_{T_{n-1}} )\,\sigma^{-\lambda} \\
&\!\!\!=\!\!\!& \sigma^\lambda\, \bigl(x(a_r|_{T_{n-1}})x^{-1}, y(b_r|_{T_{n-1}})y^{-1} \bigr)\,\sigma^{-\lambda} \\
&\!\!\!=\!\!\!& z\, (a_r|_{T_{n-1}}, b_r|_{T_{n-1}} )\,z^{-1} \\
&\!\!\!=\!\!\!& z\, (a_{r+1}|_{T_n})\,z^{-1}.
\end{eqnarray*}
Thus $w^{-1}a'_{r+1}w|_{T_n}$ is conjugate to $a_{r+1}|_{T_n}$ under~$G_n$, as desired.

\medskip
Similarly, 
for any $i\not=1,r+1$ choose an element $x\in G_{n-1}$ such that
$$u^{-1}a''_{i-1}u|_{T_{n-1}} \;=\; x\,(a_{i-1}|_{T_{n-1}})\,x^{-1}.$$
Using Proposition \ref{42G1Prop} (c) we can find an element $y\in G_{n-1}$ such that $(x,y)\in G_n$. Recall that (\ref{44SemiRigidLem1}) says that $a'_i = \sigma^\lambda\,(a''_{i-1},1)\,\sigma^{-\lambda}$ for some $\lambda\in\{0,1\}$. Thus with $z := \sigma^\lambda\,(x,y) \in G_n$ we deduce that 
\begin{eqnarray*}
w^{-1}a'_iw|_{T_n}
&\!\!\!=\!\!\!& (u,u)^{-1}\,\sigma^\lambda (a''_{i-1},1)\,\sigma^{-\lambda}\, (u,u)|_{T_n} \\
&\!\!\!=\!\!\!& \sigma^\lambda\, (u^{-1}a''_{i-1}u|_{T_{n-1}}, 1 )\,\sigma^{-\lambda} \\
&\!\!\!=\!\!\!& \sigma^\lambda\, \bigl(x(a_{i-1}|_{T_{n-1}})x^{-1}, 1 \bigr)\,\sigma^{-\lambda} \\
&\!\!\!=\!\!\!& z\, (a_{i-1}|_{T_{n-1}}, 1 )\,z^{-1} \\
&\!\!\!=\!\!\!& z\, (a_i|_{T_n})\,z^{-1}.
\end{eqnarray*}
Therefore $w^{-1}a'_iw|_{T_n}$ is conjugate to $a_i|_{T_n}$ under~$G_n$, as desired.

\medskip
Likewise, for any $1<j\le r$ choose an element $x\in G_{n-1}$ such that
$$u^{-1}b''_{j-1}u|_{T_{n-1}} \;=\; x\,(b_{j-1}|_{T_{n-1}})\,x^{-1}.$$
Using Proposition \ref{42G1Prop} (c) we can find an element $y\in G_{n-1}$ such that $(x(b_r^{-1}|_{T_{n-1}}),y)\in G_n$. Recall that (\ref{44SemiRigidLem1}) says that $b'_j = \sigma^\lambda\,(b''_{j-1},1)\,\sigma^{-\lambda}$ for some $\lambda\in\{0,1\}$.  Thus with $z := \sigma^\lambda\,(x(b_r^{-1}|_{T_{n-1}}),y) \in G_n$ we deduce that 
\begin{eqnarray*}
w^{-1}b'_jw|_{T_n}
&\!\!\!=\!\!\!& (u,u)^{-1}\,\sigma^\lambda (b''_{j-1},1)\,\sigma^{-\lambda}\, (u,u)|_{T_n} \\
&\!\!\!=\!\!\!& \sigma^\lambda\, (u^{-1}b''_{j-1}u|_{T_{n-1}}, 1 )\,\sigma^{-\lambda} \\
&\!\!\!=\!\!\!& \sigma^\lambda\, \bigl(x(b_{j-1}|_{T_{n-1}})x^{-1}, 1 \bigr)\,\sigma^{-\lambda} \\
&\!\!\!=\!\!\!& z\, ( b_rb_{j-1}b_r^{-1}|_{T_{n-1}}, 1 )\,z^{-1} \\
&\!\!\!=\!\!\!& z\, (b_j|_{T_n})\,z^{-1}.
\end{eqnarray*}
Thus $w^{-1}b'_jw|_{T_n}$ is conjugate to $b_j|_{T_n}$ under~$G_n$, as desired.

\medskip
Also, by construction we already have $w^{-1}a'_1w = (u,u)^{-1}\,\sigma\,(u,u)=\sigma=a_1.$

\medskip
Finally, for $w^{-1}b'_1w|_{T_n}$
observe that (\ref{44SemiRigidLem2}) implies that one of $c^{\pm1}|_{T_{n-1}}$ is a product of $b''_r|_{T_{n-1}}$ and some of the elements $a''_i|_{T_{n-1}}$ and $b''_j|_{T_{n-1}}$ for $i\not=r$ and $1\le j<r$ in some order.  Thus one of $u^{-1}c^{\pm1}u|_{T_{n-1}}$ is a product of $u^{-1}b''_ru|_{T_{n-1}}$ and some of the elements $u^{-1}a''_iu|_{T_{n-1}}$ and $ub''_ju^{-1}|_{T_{n-1}}$ for $i\not=r$ and $1\le j<r$ in some order. 
Here $u^{-1}b''_ru|_{T_{n-1}}$ is conjugate to $b_r|_{T_{n-1}}$ under $G_{n-1}$, so by Proposition \ref{42nHSemiDirect} (b) it lies in the coset $H_{n-1}(b_r|_{T_{n-1}})$. Also $u^{-1}a''_iu|_{T_{n-1}}$ and $ub''_ju^{-1}|_{T_{n-1}}$ for $i,j\not=r$ are conjugate to the respective $a_i|_{T_{n-1}}$ and $b_j|_{T_{n-1}}$ under $G_{n-1}$, so by the definition of $H$ they lie in~$H_{n-1}$. Together it follows that $u^{-1}cu|_{T_{n-1}}$ lies in the coset $H_{n-1}(b_r|_{T_{n-1}})$. Write $u^{-1}cu|_{T_{n-1}} = h(b_r|_{T_{n-1}})$ with $h\in H_{n-1}$. Then with $z := (h,1) \in G_n$ we deduce that 
\begin{eqnarray*}
w^{-1}b'_1w|_{T_n}
&\!\!\!=\!\!\!& (u,u)^{-1}\,(c,c^{-1})\,\sigma\, (u,u)|_{T_n} \\
&\!\!\!=\!\!\!& (u^{-1}cu|_{T_{n-1}}, u^{-1}c^{-1}u|_{T_{n-1}} )\,\sigma \\
&\!\!\!=\!\!\!& \bigl( h(b_r|_{T_{n-1}}), (b_r|_{T_{n-1}})^{-1}h^{-1} \bigr)\,\sigma \\
&\!\!\!=\!\!\!& z\, (b_r|_{T_{n-1}}, b^{-1}_r|_{T_{n-1}} )\,z^{-1} \\
&\!\!\!=\!\!\!& z\, (b_1|_{T_n})\,z^{-1}.
\end{eqnarray*}
Thus $w^{-1}b'_1w|_{T_n}$ is conjugate to $b_1|_{T_n}$ under~$G_n$, as desired.

\medskip
Together this shows that ($*_{n-1}$) implies ($*_n$). Thus by induction ($*_n$) is true for all $n\ge0$, finishing the proof of Theorem \ref{44SemiRigid}.
\end{Proof}

\begin{Prop}\label{44WAiBiConj1}
For any $1\le i\le r$, the elements $a_i$ and $b_i$ are conjugate under~$W$, but not under~$G$.
\end{Prop}

\begin{Proof}
The recursion relations (\ref{4RecRels}) directly show that $b_1=(b_r,b_r^{-1})\,\sigma\sim\sigma=a_1$, and if $1<i\le r$ with $b_{i-1}\sim a_{i-1}$, they show that $b_i = (b_rb_{i-1}b_r^{-1},1) \sim (b_{i-1},1) \sim (a_{i-1},1) =a_i$. By induction we deduce that $b_i\sim a_i$ for all $1\le i\le r$.

Suppose that for some $1\le i\le r$ the elements $a_i$ and $b_i$ are conjugate under~$G$. Let $i$ be minimal with this property and choose $z\in G$ with $b_i=za_iz^{-1}$. By Proposition \ref{42HSemiDirect} (c) we can write $z=(h,h')\,(b_r,b_r)^\lambda\,\sigma^\mu$ with $h,h'\in H$ and $\lambda,\mu\in\{0,1\}$. If $i=1$, we then have
$$(b_r,b_r^{-1})\,\sigma = b_1 = za_1z^{-1} = z\sigma z^{-1} = (h,h')\,\sigma\,(h,h')^{-1} = (hh^{\prime-1},h'h^{-1})\,\sigma$$
and therefore $b_r = hh^{\prime-1} \in H$. But this contradicts Proposition \ref{42HSemiDirect} (b); hence $i>1$. If $\mu=1$, the conjugate $za_iz^{-1} = z\,(a_{i-1},1)\,z^{-1}$ has the form $(1,*)$; since $b_{i-1}\not=1$, it is therefore different from $b_i=(b_rb_{i-1}b_r^{-1},1)$. Thus $\mu=0$, and hence
$$(b_rb_{i-1}b_r^{-1},1) = b_i = za_iz^{-1} = 
(hb_r^\lambda,h'b_r^\lambda)\,(a_{i-1},1)\,(hb_r^\lambda,h'b_r^\lambda)^{-1}
= (hb_r^\lambda a_{i-1}(hb_r^\lambda)^{-1},1).$$
This shows that $b_{i-1} = (b_r^{-1}hb_r^\lambda) a_{i-1}(b_r^{-1}hb_r^\lambda)^{-1}$, and so $b_{i-1}$ is conjugate to $a_{i-1}$ under~$G$. By the minimality of $i$ this yields a contradiction, proving that $a_i$ and $b_i$ are not conjugate under~$G$ for any $1\le i\le r$.
\end{Proof}

\begin{Prop}\label{44WAiBiConj2}
For any $1\le i\le r$, any element of $G$ which is conjugate to $a_i$ or $b_i$ under $W$ is conjugate to precisely one of $a_i$ and $b_i$ under~$G$.
\end{Prop}

\begin{Proof}
In view of Proposition \ref{44WAiBiConj1} it remains to prove that any element $z\in G$ which is conjugate to $a_i$ under $W$ is conjugate to at least one of $a_i$ and $b_i$ over~$G$. Again we will show this by induction over~$i$.

If $i=1$, we have $a_1=\sigma$ and hence $z=(x,x^{-1})\,\sigma$ for some $x\in W$. 
Since $z\in G$, Proposition \ref{42G1Prop} (a) implies that $x\in G$. By Proposition \ref{42HSemiDirect} (b) we can thus write $x=h$ or $x=hb_r$ with $h\in H$. Then $w := (h,1)$ lies in $G$ by Proposition \ref{42HSemiDirect} (c), and according to the case we deduce that
\begin{eqnarray*}
z &\!\!=\!\!& (h,h^{-1})\,\sigma \, =\, (h,1)\,\sigma\,(h,1)^{-1} 
\,=\, wa_1w^{-1} \qquad\hbox{or}\\
z &\!\!=\!\!& (hb_r,b_r^{-1}h^{-1})\,\sigma 
\,=\, (h,1)\,(b_r,b_r^{-1})\,\sigma\,(h,1)^{-1} \,=\, wb_1w^{-1}.
\end{eqnarray*}
Thus $z$ is conjugate to $a_1$ or $b_1$ under~$G$, as desired.

If $1<i\le r$, after possibly replacing $z$ by $a_1za_1^{-1}$, we may assume that $z$ is conjugate to $a_i$ under the subgroup $W\times W \subset W$. Since $a_i = (a_{i-1},1)$, this means that $z=(x,1)$ where $x\in W$ is conjugate to $a_{i-1}$ under~$W$. As $z\in G$, Proposition \ref{42G1Prop} (a) implies that $x\in G$. By the induction hypothesis we thus know that $x$ is conjugate to $a_{i-1}$ or $b_{i-1}$ under~$G$. Choose an element $u\in G$ with $x=ua_{i-1}u^{-1}$ or $x=ub_rb_{i-1}b_r^{-1}u^{-1}$. By Proposition \ref{42G1Prop} (c) there exists an element $v\in G$ such that $w := (u,v)\in G$.
According to the case we deduce that
\begin{eqnarray*}
z &\!\!=\!\!& (ua_{i-1}u^{-1},1) \, =\, (u,v)\,(a_{i-1},1)\,(u,v)^{-1} 
\,=\, wa_iw^{-1} \qquad\hbox{or}\\
z &\!\!=\!\!& (ub_rb_{i-1}b_r^{-1}u^{-1},1)
\,=\, (u,v)\,(b_rb_{i-1}b_r^{-1},1)\,(u,v)^{-1} \,=\, wb_iw^{-1}.
\end{eqnarray*}
Thus $z$ is conjugate to $a_i$ or $b_i$ under~$G$, as desired.
\end{Proof}

\begin{Rem}\label{44WAiBiConj3}
\rm Proposition \ref{44WAiBiConj2} does not directly extend to $i>r$, for instance because $(b_r,b_r)\in G$ is conjugate to $a_{r+1}=(a_r,b_r)$ under $W$ but not under~$G$, as can easily be shown.
\end{Rem}


\subsection{Normalizer}
\label{45Normalizer}

Next we will determine the normalizer
\UseTheoremCounterForNextEquation
\begin{equation}\label{45NormDef}
N(r) := \Norm_W(G(r)).
\end{equation}
which we also abbreviate by $N$ until the end of this subsection.

\begin{Lem}\label{45Lem1}
The group $N$ normalizes $H$.
\end{Lem}

\begin{Proof}
Direct consequence of Proposition \ref{42HSemiDirect} (a) and the fact that $\sgn_r$ is defined on all of~$W$.
\end{Proof}

\medskip
Let $\diag\colon W\to W\times W$, $w\mapsto (w,w)$ denote the diagonal embedding.

\begin{Lem}\label{45Lem2}
\begin{enumerate}
\item[(a)] We have $\diag(G) \subset G$.
\item[(b)] We have $N = G\cdot\langle(1,b_r)\rangle\cdot\diag(N)$.
\end{enumerate}
\end{Lem}

\begin{Proof}
(a) is a direct consequence of Proposition \ref{42HSemiDirect} (b) and (c).
For (b) consider an arbitrary element $(u,v)\in W\times W$. If $(u,v)$ lies in~$N$, it normalizes $G^1 = G\cap(W\times W)$. Then by Proposition \ref{42G1Prop} (c) both $u$ and $v$ normalize $\proj_1(G^1)=\proj_2(G^1)= \nobreak G$, in other words we have $u,v\in N$.

Conversely assume that $u,v\in N$. Then $(u,v)$ already normalizes $H\times H$ by Lemma \ref{45Lem1}. By Proposition \ref{42HSemiDirect} (c) it therefore normalizes $G$ if and only if both 
\begin{eqnarray*}
(u,v)\,(b_r,b_r)\,(u,v)^{-1} &\!\!=\!\!& (ub_ru^{-1},vb_rv^{-1}) \qquad\hbox{and} \\
(u,v)\,\sigma\,(u,v)^{-1} &\!\!=\!\!& (uv^{-1},vu^{-1})\,\sigma
\end{eqnarray*}
lie in~$G$. 
Since $u$ and $v$ already normalize $G$ and~$H$, they normalize the complement $G\setminus H = Hb_r$; hence both $ub_ru^{-1}$ and $vb_rv^{-1}$ lie in~$Hb_r$. Thus the first of the above two elements automatically lies in $Hb_r\times Hb_r \subset G$. 
By Proposition \ref{42HSemiDirect} (c) the second lies in $G$ if and only if $(uv^{-1},vu^{-1}) \in Hb_r^\lambda\times Hb_r^\lambda$ for some $\lambda\in\{0,1\}$. But if $vu^{-1} \in Hb_r^\lambda$, it already follows that $uv^{-1} = (vu^{-1})^{-1}\in b_r^{-\lambda}H = Hb_r^\lambda$. Thus the second element lies in $G$ if and only if $vu^{-1} \in Hb_r^\lambda$ for some $\lambda\in\{0,1\}$.

Together we find that an element $(u,v)\in W\times W$ lies in $N$ if and only if $u\in N$ and $v\in H\langle b_r\rangle u$. Since $\sigma\in G$ and $(1\times H)\subset G$ by Proposition \ref{42HSemiDirect} (c), we deduce that 
\begin{eqnarray*}
N &\!\!=\!\!& G\cdot(N\cap(W\times W)) \\
&\!\!=\!\!& G\cdot(1\times H\langle b_r\rangle)\cdot\diag(N) \\
&\!\!=\!\!& G\cdot(1\times H) \cdot\langle (1,b_r)\rangle\cdot\diag(N) \\
&\!\!=\!\!& G\cdot\langle (1,b_r)\rangle\cdot\diag(N)
\end{eqnarray*}
as desired.
\end{Proof}

\medskip
Now we recursively define elements 
\UseTheoremCounterForNextEquation
\begin{equation}\label{45WiDef}
\biggl\{\!\begin{array}{lll}
w_1     &\!\!\!:=\, (1,b_r), & \\[3pt]
w_{i+1} &\!\!\!:=\, (w_i,w_i) & \hbox{for all $i\ge1$.}\\
\end{array}\!\biggr\}
\end{equation}
By induction Lemma \ref{45Lem2} (b) implies that $w_i\in N$ for all $i\ge1$. Also, since $b_r$ has order two, by induction the same follows for all~$w_i$. Moreover, by induction we find that the restriction $w_i|_{T_n}$ is trivial for all $i\ge n$; hence the sequence $w_1,w_2,\ldots$ converges to $1$ within~$N$. Thus the following map is well-defined:
\UseTheoremCounterForNextEquation
\begin{equation}\label{45MapDef}
\textstyle \phi\colon
\prod\limits_{i=1}^\infty\BF_2 \longto N, \quad (k_1,k_2,\ldots) \mapsto w_1^{k_1}w_2^{k_2}\cdots.
\end{equation}
By construction it is continuous, but not a homomorphism. Note that it satisfies the basic formula
\UseTheoremCounterForNextEquation
\begin{equation}\label{45MapForm}
\phi(k_1,k_2,\ldots) \ =\ w_1^{k_1}\cdot\diag(\phi(k_2,k_3,\ldots)).
\end{equation}

\begin{Lem}\label{45Lem3}
The map $\phi$ induces a homomorphism $\bar\phi\colon \prod_{i=1}^\infty\BF_2 \to N/G$.
\end{Lem}

\begin{Proof}
It suffices to show that for all $i>j\ge1$ the images of $w_i$ and $w_j$ in $N/G$ commute with each other; in other words that the commutator $[w_i,w_j]$ lies in~$G$. We will prove this by induction on~$j$.

For $j=1$ we have $[w_i,w_1] = [(w_{i-1},w_{i-1}),(1,b_r)] = (1,[w_{i-1},b_r])$. Here $[w_{i-1},b_r]$ lies in~$G$, because $w_{i-1}\in N$. Being a commutator, this element also lies in the kernel of $\sgn_r$; hence by Proposition \ref{42HSemiDirect} (a) it lies in~$H$. Thus $[w_i,w_1]$ lies in $1\times H$ and hence in $G$ by Proposition \ref{42HSemiDirect} (c), as desired.

For $j>1$ we have $[w_i,w_j] = \diag([w_{i-1},w_{j-1}])$. Here $[w_{i-1},w_{j-1}]$ lies in~$G$ by the induction hypothesis. Thus $[w_i,w_j]$ lies in $G$ by Lemma \ref{45Lem2} (a), as desired.
\end{Proof}

\begin{Lem}\label{45Lem4}
The homomorphism $\bar\phi$ is injective.
\end{Lem}

\begin{Proof}
Let $(k_1,k_2,\ldots)$ be an element of the kernel of~$\bar\phi$. Then $w := \phi(k_1,k_2,\ldots)$ lies in~$G$. By the formula (\ref{45MapForm}) we have $w=w_1^{k_1}\cdot \diag(u) = (u,b_r^{k_1}u)$ for $u := \phi(k_2,k_3,\ldots)$, and so this element lies already in~$G^1$. By Proposition \ref{42HSemiDirect} (d) this requires that $u\in G$ and $b_r^{k_1}\in H$. By Proposition \ref{42HSemiDirect} (b) we must therefore have $k_1=0$. Moreover, the fact that $u\in G$ means that $(k_2,k_3,\ldots)$ also lies in the kernel of~$\bar\phi$.

For every element $(k_1,k_2,\ldots)$ of the kernel of~$\bar\phi$ we have thus proved that $k_1=0$ and that $(k_2,k_3,\ldots)$ again lies in the kernel of~$\bar\phi$. By an induction on $i$ we can deduce from this that for every $i\ge1$ and every element $(k_1,k_2,\ldots)$ of the kernel of~$\bar\phi$ we have $k_i=0$. This means that the kernel of $\bar\phi$ is trivial, and so the homomorphism $\bar\phi$ is injective, as desired.
\end{Proof}

\begin{Lem}\label{45Lem5}
The homomorphism $\bar\phi$ is surjective.
\end{Lem}

\begin{Proof}
The assertion is equivalent to $N = G\cdot\phi\bigl(\prod_{i=1}^\infty\BF_2\bigr)$. For this it suffices to show that $N_n = G_n\cdot\phi\bigl(\prod_{i=1}^\infty\BF_2\bigr)|_{T_n}$ for all~$n\ge0$. This is trivial for $n=0$, so assume that $n>0$ and that the equality holds for~$n-1$. 
Using, in turn, Lemma \ref{45Lem2} (b), the induction hypothesis, Lemma \ref{45Lem2} (a), and the formula (\ref{45MapForm}) we deduce that
\begin{eqnarray*}
N_n &=& G_n\cdot\langle(1,b_r)|_{T_n}\rangle\cdot\diag(N_{n-1}) \\
&=& \textstyle G_n\cdot\langle w_1|_{T_n}\rangle\cdot\diag\bigl( G_{n-1} \cdot 
\phi\bigl(\prod_{i=1}^\infty\BF_2\bigr)|_{T_{n-1}} \bigr) \\
&=& \textstyle G_n\cdot\langle w_1|_{T_n}\rangle\cdot\diag(G_{n-1}) \cdot 
\diag\bigl(\phi\bigl(\prod_{i=1}^\infty\BF_2\bigr)|_{T_{n-1}} \bigr) \\
&=& \textstyle G_n\cdot\langle w_1|_{T_n}\rangle \cdot 
\diag\bigl(\phi\bigl(\prod_{i=1}^\infty\BF_2\bigr)|_{T_{n-1}} \bigr) \\
&=& \textstyle G_n\cdot \phi\bigl(\prod_{i=1}^\infty\BF_2\bigr)|_{T_n},
\end{eqnarray*}
so the equality holds for~$n$. Thus it follows for all $n\ge0$ by induction, and we are done.
\end{Proof}

\medskip
Combining Lemmas \ref{45Lem3} through \ref{45Lem5} now implies:

\begin{Thm}\label{45NormThm}
The map $\phi$ induces an isomorphism 
$$\textstyle \bar\phi\colon \prod\limits_{i=1}^\infty\BF_2 \stackrel{\sim}{\longto} N/G.$$
\end{Thm}


\subsection{Inclusions}
\label{47Inclusions}

In this subsection we determine the possible inclusions between the groups $G(r)$ and $N(r)$ studied above for different values of~$r$.

\begin{Thm}\label{47IncluThmAll}
For any $r>r'\ge1$ we have $G(r')\subset N(r') \subset G(r)$.
\end{Thm}

\begin{Rem}\label{47NoIncluRem}
\rm The formula $1-2^{-r}$ for the Hausdorff dimension of~$G(r)$, which is strictly monotone increasing with $r$, shows that for $r'<r$ not even a conjugate of $G(r)$ can be contained in~$G(r')$. The same holds for $N(r')$ in place of $G(r')$, because with \ref{45NormThm} one can easily show that their Hausdorff dimensions are equal.
\end{Rem}

\begin{Proof}
By induction it suffices to prove Theorem \ref{47NoIncluRem} when $r'=r-1$. To distinguish the respective generators, we endow the generators of $G(r')$ with a prime~$'$.
They are thus the elements of $W$ that are uniquely determined by the equations
\UseTheoremCounterForNextEquation
\begin{equation}\label{47RecRelsPrime}
\left\{\begin{array}{lcll}
a'_1     \!\!\!&=&\!\!\! \sigma & \\[3pt]
a'_r \!\!\!&=&\!\!\! (a'_{r-1},b'_{r-1}) & \\[3pt]
a'_i     \!\!\!&=&\!\!\! (a'_{i-1},1) & \hbox{\ for all $i>1$ with $i\not=r$,}
 \\[3pt]
b'_1     \!\!\!&=&\!\!\! (b'_r,b_r^{\prime-1})\,\sigma & \\[3pt]
b'_j     \!\!\!&=&\!\!\! (b'_rb'_{j-1}b_r^{\prime-1},1) & \hbox{\ for all $1<j\le r-1$.} \\[3pt]
\end{array}\right\}
\end{equation}
Moreover, by Theorem \ref{45NormThm} for $r-1$ in place of~$r$ the group $N(r')$ is topologically generated by these together with the elements recursively defined by 
\UseTheoremCounterForNextEquation
\begin{equation}\label{45WiDefPrime}
\biggl\{\!\begin{array}{lll}
w'_1     &\!\!\!:=\, (1,b'_{r-1}), & \\[3pt]
w'_{i+1} &\!\!\!:=\, (w'_i,w'_i) & \hbox{for all $i\ge1$.}\\
\end{array}\!\biggr\}
\end{equation}
We must prove that all these elements lie in $G(r)$. For this we do not need to mention the group $G(r')$ at all, and can use the previously established results about $G(r)$.
For simplicity we again abbreviate $G :=G(r)$. Note that by Proposition \ref{41GenSigns} with $r-1$ in place of $r$ we have $\sgn_r(a'_i) = \sgn_r(b'_j)=1$ for all $i$ and~$j$.

\begin{Lem}\label{47IncluLem1}
For all $1\le j\le r-1$ we have $b'_j\in G(r)$.
\end{Lem}

\begin{Proof}
It suffices to prove that for all $n\ge0$ we have $b'_j|_{T_n}\in G_n$ for all $j$.
For $n\le r$ this follows from the fact that $G_n=W_n$ by Proposition \ref{41GnAll}. So assume that it is true for some $n\ge r$ and all~$j$. Then $\sgn_r$ factors through~$W_n$, and since $\sgn_r(b'_j)=1$ for all~$j$, we deduce that $b'_j|_{T_n}\in H_n$ for all $j$. The recursion relations (\ref{47RecRelsPrime}) thus imply that all $b'_j|_{T_{n+1}}$ lie in $(H_n\times H_n)\rtimes\langle\sigma\rangle$ and hence in $G_{n+1}$ by Proposition \ref{42nHSemiDirect} (c). By induction the desired assertion therefore holds for all $n\ge0$.
\end{Proof}

\begin{Lem}\label{47IncluLem2}
For all $i\ge1$ we have $a'_i\in G$.
\end{Lem}

\begin{Proof}
Since the recursion relations (\ref{47RecRelsPrime}) for $a'_1,\ldots,a'_{r-1}$ coincide with the recursion relations (\ref{4RecRels}) for $a_1,\ldots,a_{r-1}$, we have $a'_i=a_i\in G$ for all $1\le i\le r-1$. In particular we have $a'_{r-1}\in G$ and $b'_{r-1}\in G$ by Lemma \ref{47IncluLem1}. Since $\sgn_r(a'_{r-1}) = \sgn_r(b'_{r-1})=1$, they are actually contained in $H$, and so $a'_r = (a'_{r-1},b'_{r-1})$ is contained in $H\times H\subset G$.
Now suppose we know that $a'_i\in G$ for some $i\ge r$. Then the fact that $\sgn_r(a_i')=1$ implies that actually $a'_i\in H$. Thus $a_{i+1} = (a_i,1) \in H\times 1\subset G$. By induction it follows that $a'_i\in G$ for all $i\ge r$, and we are done.
\end{Proof}

\begin{Lem}\label{47IncluLem3}
For all $i\ge1$ we have $w'_i\in G$.
\end{Lem}

\begin{Proof}
By Lemma \ref{47IncluLem1} we have $b'_{r-1}\in G$. Since moreover ${\sgn_r(b'_{r-1})=1}$, this element already lies in~$H$. Thus $w'_1=(1,b'_{r-1})$ lies in $1\times H$ and hence in~$G$. Also, if $w'_i\in G$ for some $i\ge1$, then $w'_{i+1} = \diag(w'_i) \in G$ by Lemma \ref{45Lem2} (a). By induction it follows that $w'_i\in G$ for all $i\ge1$.
\end{Proof}

\medskip
By combining Lemmas \ref{47IncluLem1} through \ref{47IncluLem3} we deduce Theorem \ref{47IncluThmAll}.
\end{Proof}


\subsection{Iterated monodromy groups}
\label{48Monodromy}

Now let $f$ be a quadratic morphism over a field $k$ of characteristic $\not=2$ with the postcritical orbit $P\subset\BP^1(\bar k)$. (Compare the introduction and \cite[\S\ref{17Monodromy}]{Pink2013b}.) In this subsection we assume that $P$ is infinite.
Let $\rho\colon {\pi_1^\et(\BP^1_k \setminus P)} \to W$ be the monodromy representation and $G^\geom\subset G^\arith\subset W$ the geometric and arithmetic fundamental groups associated to~$f$. Let $p_0$, $q_0\in\BP^1(\bar k)$ denote the two critical points of~$f$. Note that necessarily $f(p_0)\not=f(q_0)$. 

\begin{Thm}\label{48GgeomThm}
\begin{enumerate}
\item[(a)] If $f^{r+1}(p_0)\not=f^{r+1}(q_0)$ for all $r\ge1$, then $G^\geom=G^\arith=W$.
\item[(b)] Otherwise let $r\ge1$ be minimal with $f^{r+1}(p_0)=f^{r+1}(q_0)$. Then there exists $w\in W$ such that $G^\geom=wG(r)w^{-1}$.
\end{enumerate}
\end{Thm}

\begin{Proof}
Set $X:=\BP^1(\bar k)$ with the induced map $f\colon X\to X$, the set of critical points $C := \{p_0,q_0\}$, and the postcritical orbit $P := \bigcup_{n\ge1}f^n(C)$. By the recursion relations (\ref{0GenRecConj1}) and the product relation, the given generators $b_p$ for $p\in P$ are strongly of type $(X,f,C)$ in the sense of Subsection \ref{41TheGroups}. 
In the case (a), Propositions \ref{41LargeAB} and \ref{41LargeC} imply that $G^\geom=W$. Since in general $G^\geom\subset G^\arith\subset W$, these inclusions are equalities in this case, proving (a).

In the case (b) the generators satisfy the recursion relations (\ref{4GenRecConj2}), which up to renaming the generators are the same those in Theorem \ref{44SemiRigid}. Thus Theorem \ref{44SemiRigid} shows that $G^\geom=wG(r)w^{-1}$ for some $w\in W$.
\end{Proof}

\medskip
In the rest of this subsection we consider only the case \ref{48GgeomThm} (b). For simplicity we change the identification of trees used in \cite[\S\ref{17Monodromy}]{Pink2013b} by the automorphism~$w$, after which we have $G^\geom=G(r)$, and the given generators $b_{p_i}$ and $b_{q_j}$ of $G^\geom$ are conjugate to $a_i$ and $b_j$ under $G(r)$ for all $i\ge1$ and $1\le j\le r$. Then $G^\arith$ is contained in the normalizer $N(r)$ of $G(r)$, and to describe it it suffices to describe the factor group $G^\arith/G(r) \subset N(r)/G(r)$. More precisely we will determine the composite homomorphism
\UseTheoremCounterForNextEquation
\begin{equation}\label{48GalRep}
\bar\rho\colon \Gal(\bar k/k) \onto G^\arith/G(r) \into N(r)/G(r)
\end{equation}
obtained from the homomorphism~$\rho$.
Its composite with the isomorphism from Theorem \ref{45NormThm} must have the form
\UseTheoremCounterForNextEquation
\begin{equation}\label{48CharDecomp}
(\chi_1,\chi_2,\ldots) :\ \Gal(\bar k/k) \stackrel{\bar\rho}{\longto} N(r)/G(r) \,\cong\, \prod_{i=1}^\infty\BF_2
\end{equation}
for continuous homomorphisms $\chi_i\colon \Gal(\bar k/k) \to \nobreak \BF_2$ that remain to be determined.

\begin{Lem}\label{48LemDiag}
The homomorphisms $\chi_i$ for $i\ge1$ are all equal.
\end{Lem}

\begin{Proof}
Combining the short exact sequence (\ref{17ShortExact}) and the diagrams (\ref{17RecArithDiag}) and (\ref{17RecGeomDiag}) of \cite{Pink2013b}, we find the left third of the following commutative diagram:
\UseTheoremCounterForNextEquation
\begin{equation}\label{17RecFactDiag}
\vcenter{\xymatrix@R+0pt@C+10pt{
\Gal(\bar k/k) \ar@{->>}[r]^-{\bar\rho} &
G^\arith/G(r) \ar@{^{ (}->}[r] & N(r)/G(r) & 
\prod\limits_{i=1}^\infty\!\BF_2 \ar[dd]_{(*)} \ar[l]^-\sim_-{\bar\phi} \\
\Gal(\bar k/k) \ar@{=}[u] \ar[r] \ar@{=}[d] &
{\frac{\textstyle G^\arith\cap(W\times W)}{\textstyle G(r)\cap(W\times W)}}
\ar[u]^-\wr \ar[d]_-{[\proj_1]} \ar@{^{ (}->}[r] & 
{\frac{\textstyle N(r)\cap(W\times W)}{\textstyle G(r)\cap(W\times W)}}
\ar[u]^-\wr \ar[d]_-{[\proj_1]} 
& \\
\Gal(\bar k/k) \ar@{->>}[r]^-{\bar\rho} &
G^\arith/G(r) \ar@{^{ (}->}[r] & N(r)/G(r) & 
\prod\limits_{i=1}^\infty\!\BF_2 \ar[l]^-\sim_-{\bar\phi} \\}}
\end{equation}
The middle third comes from the inclusions $G(r)\subset G^\arith\subset N(r)$. The homomorphisms from the second row to the first are obtained from the natural inclusions, and they are isomorphisms, because $G(r)$ acts transitively on level~$1$ and therefore $W = G(r)\cdot(W\times W)$. In the right third the isomorphism $\bar\phi$ is that from Theorem \ref{45NormThm}.
The rightmost vertical homomorphism marked $(*)$ is defined to make everything commute.
To determine it recall that $(k_1,k_2,\ldots)\in\prod_{i=1}^\infty\!\BF_2$ corresponds to the element $\phi(k_1,k_2,\ldots)\in N$, which by the definition (\ref{45WiDef}) of $w_1$ and the formula (\ref{45MapForm}) is equal to $(1,b_r)^{k_1}\cdot\diag(\phi(k_2,k_3,\ldots))$.
Thus it lies in $W\times W$, and its image under $\proj_1$ is $\proj_1\bigl(\diag(\phi(k_2,k_3,\ldots))\bigr) = \phi(k_2,k_3,\ldots)$. It follows that to make the diagram (\ref{17RecFactDiag}) commutes the map $(*)$ must be defined by $(k_1,k_2,\ldots) \mapsto (k_2,k_3,\ldots)$. This implies that $\chi_1=\chi_2=\chi_3=\ldots$, as desired.
\end{Proof}

\begin{Lem}\label{48LemChi1}
Identify $\BF_2$ with the the symmetric group~$S_2$ on two letters, and let $\theta\colon \allowbreak\Gal(\bar k/k) \to S_2\cong\BF_2$ denote the homomorphism describing the action of $\Gal(\bar k/k)$ on the set of critical points $C = \{p_0,q_0\} \subset \BP^1(\bar k)$.
Then $\chi_1=\theta$.
\end{Lem}

\begin{Proof}
By Proposition \ref{44WAiBiConj1} the elements $a_1$ and $b_1$ lie in two different conjugacy classes of~$G(r)$. Since $w_1\in N(r)$ is an element of order two satisfying $w_1^{-1}a_1w_1=(1,b_r^{-1})\,\sigma\,(1,b_r) = (b_r,b_r^{-1})\,\sigma = b_1$, conjugation by $w_1$ interchanges these two conjugacy classes.
Also, for any $i\ge1$ the element $w_{i+1}=(w_i,w_i)$ commutes with ${a_1=\sigma}$, 
and since by Lemma \ref{45Lem3} it commutes with $w_1$ modulo~$G(r)$, it follows that conjugation by $w_{i+1}$ maps each of the conjugacy classes of $a_1$ and $b_1$ to themselves.
Together this implies that conjugation by $\tilde w := w_1w_2w_3\cdots = \phi(1,1,1,\ldots)$ interchanges the conjugacy classes of $a_1$ and $b_1$ in~$G(r)$.

Let us name the generators of $G^\geom$ as in Theorem \ref{44SemiRigid}. Then $a'_1$ and $b'_1$ are generators of the images of inertia groups at $p_1$ and~$q_1$, respectively. 
{}From Theorem \ref{44SemiRigid} we know that $a'_1$ is conjugate to $a_1$ under~$G(r)$, and $b'_1$ is conjugate to $b_1$ under~$G(r)$. Thus conjugation by $\tilde w$ interchanges the two distinct conjugacy classes of $a'_1$ and $b'_1$ in~$G(r)$.

Now consider any element $\gamma\in \pi_1^\et(\BP^1_k \setminus P)$. By Lemma \ref{48LemDiag} its image in $N(r)$ has the form $\rho(\gamma) = g\tilde w^\lambda$ for some $g\in G(r)$ and $\lambda\in\{0,1\}$. Then $\lambda=1$ if and only if conjugation by $\rho(\gamma)$ interchanges the conjugacy classes of $a'_1$ and $b'_1$, which it does if and only if the Galois action of $\gamma$ interchanges the critical points $p_0$ and~$q_0$. This implies that $\chi_1=\theta$, as desired.
\end{Proof}

\medskip
By combining Lemmas \ref{48LemDiag} and \ref{48LemChi1} we deduce:

\begin{Thm}\label{48GarithThm1}
In the case \ref{48GgeomThm} (b), the homomorphism $\bar\rho\colon \Gal(\bar k/k) \to N(r)/G(r)\cong \prod_{i=1}^\infty\BF_2$ is equal to the homomorphism $\theta\colon \Gal(\bar k/k)\to S_2\cong \BF_2$ describing the action on the critical points $C = \{p_0,q_0\} \subset \BP^1(\bar k)$ followed by the diagonal embedding $\BF_2\into \prod_{i=1}^\infty\BF_2$.
\end{Thm}

Now consider the element $\tilde w := \phi(1,1,1,\ldots)\in N(r)$, which by (\ref{45WiDef}) and (\ref{45MapForm}) satisfies
\UseTheoremCounterForNextEquation
\begin{equation}\label{48GTilde}
\tilde w = w_1\cdot\diag(\tilde w)
= (1,b_r)\,(\tilde w,\tilde w)
= (\tilde w,b_r\tilde w).
\end{equation}
Note that this recursion relation determines $\tilde w$ in its own right by \cite[Prop.\ \ref{15RecRelsProp}]{Pink2013b}. Let $\tilde G(r)$ denote the subgroup of $N(r)$ that is generated by $G(r)$ and~$\tilde w$. As the image of $\tilde w$ in $N(r)/G(r)$ has order $2$, we have $[\tilde G(r):G(r)]=2$. 
Theorem \ref{48GarithThm1} now directly implies:


\begin{Thm}\label{48GarithThm2}
In the case \ref{48GgeomThm} (b) we have $G^\arith\sim G(r)$ if the two critical points of $f$ are $k$-rational, and $G^\arith\sim\tilde G(r)$ otherwise.
\end{Thm}

\medskip
Consider now in addition an intermediate field $k\subset k'\subset\bar k$ and a point $x'\in \BP^1(k')\setminus P$. Viewing $x'$ as a point in $\BP^1(\bar k)$, we can construct the associated regular rooted binary tree $T_{x'}$ with set of vertices $\coprod_{n\ge0}f^{-n}(x')$, as in \cite[\S\ref{17Monodromy}]{Pink2013b} with $\bar k$ in place of~$L$. This time we are interested in the natural action of $\Gal(\bar k/k')$ on~$T_{x'}$. Identify $T_{x'}$ with the standard tree $T$ in some way and let $G_{x'} \subset W$ denote the image of the continuous homomorphism $\Gal(\bar k/k')\to W$ describing the action on~$T_{x'}$. 

\begin{Thm}\label{47SpecDecomp}
In the case \ref{48GgeomThm} (b) there exists $w\in W$ such that 
$$G_{x'}\subset w\tilde G(r)w^{-1}.$$
\end{Thm}

\begin{Proof}
By \cite[Prop.\ \ref{18SpecDecomp}]{Pink2013b} there exists $w\in W$ such that $G_{x'}\subset wG^\arith w^{-1}$, so the theorem results from Theorem \ref{48GarithThm2}.
\end{Proof}

\medskip
This partially answers a question of Rafe Jones \cite[Question 3.4]{JonesSurvey2013}.


\subsection{Specialization}
\label{49MonodromyIncl}

Let $f$ be a quadratic morphism over a field $k$ of characteristic $\not=2$, with the critical points $p_0$, $q_0\in\BP^1(\bar k)$. In this subsection we assume only that \UseTheoremCounterForNextEquation
\begin{equation}\label{49PostCritRel}
f^{r+1}(p_0) = f^{r+1}(q_0)
\end{equation}
for some $r\ge1$. Let $G^\geom\subset G^\arith\subset W$ be the geometric and arithmetic fundamental groups associated to~$f$. 

As before let $P\subset\BP^1(\bar k)$ denote the postcritical orbit of~$f$. If $P$ is infinite, the preceding section determines $G^\geom$ and $G^\arith$ and gives an upper bound for $G_{x'}$ for any point $x'\in \BP^1(k')\setminus P$, up to conjugacy under~$W$.
If $P$ is finite, the necessary knowledge on $G^\geom$ is not yet available to deduce the same kind of consequences for $G^\arith$ and~$G_{x'}$. Nevertheless, using specialization we can at least find similar upper bounds as in the infinite case. 

As a preparation we show:

\begin{Lem}\label{49pSpecLem}
There exist a noetherian normal integral domain $R$ whose residue field $k'$ is a finite extension of~$k$, and a morphism $F: \BP^1_{\Spec R} \to \BP^1_{\Spec R}$ which is fiberwise of degree~$2$, such that the quadratic morphism over $k'$ induced by $F$ coincides with the base change of~$f$, and the quadratic morphism induced by $F$ over the quotient field of $R$ has an infinite postcritical orbit satisfying the analogue of (\ref{49PostCritRel}).
\end{Lem}

\begin{Proof}
Choose $\underline{t}=(t_1,\ldots,t_6)\in k^6$ such that $f(x) = {(t_1x^2+t_2x+t_3)/}{(t_4x^2+t_5x+t_6)}$. Let $T_1,\ldots,T_6$ denote the variables on the affine space $\BA^6_k$. Then the points in $\BA^6_k$ where $F(x) := {(T_1x^2+T_2x+T_3)/}{(T_4x^2+T_5x+T_6)}$ defines a quadratic morphism form a Zariski open neighborhood $U\subset \BA^6_k$ of~$\underline{t}$.
Let $X\subset \BP_1\times U$ denote the closed subscheme defined by the equation $dF=0$. Thus $X$ consists of the critical points of $F$ and is therefore a finite \'etale Galois covering of degree~$2$ of~$U$. In particular it is smooth of dimension $6$ over~$k$.
Let $P_0:X\into \BP^1\times X$ denote the section coming from the tautological embedding $X\into \BP^1\times U$. Let $Q_0:X\into \BP^1\times X$ denote the section  obtained by twisting $P_0$ with the nontrivial Galois automorphism of $X$ over~$U$. Then the equation $F^{r+1}(P_0) = F^{r+1}(Q_0)$ makes sense everywhere over~$X$. Being a non-trivial equation on the irreducible smooth scheme $\BP^1\times X$, the points satisfying it form a closed subscheme $Y\subset X$ which is equidimensional of dimension~$5$. 

By the assumption (\ref{49PostCritRel}) there exists a point $y\in Y$ over $\underline{t}\in U$. Choose any irreducible component $Y'$ of $Y$ which contains~$y$. Let $\tilde Y'$ be the normalization of $Y'$ and $\tilde y\in \tilde Y'$ a point over~$y$.
Let $R$ be the local ring of $\tilde Y'$ at~$\tilde y$, which by construction is a noetherian normal integral domain whose residue field $k'$ is a finite extension of~$k$. For simplicity denote the quadratic morphism over $\Spec R$ induced by~$F$ again by~$F$. Then the quadratic morphism over $k'$ induced by $F$ coincides with that induced by~$f$. Moreover, the quotient field $K$ of $R$ is the function field of the irreducible variety $Y'$ over~$k$. By construction the quadratic morphism $F$ over $K$ satisfies the analogue of (\ref{49PostCritRel}).

It remains to show that the postcritical orbit of $F$ over $K$ is infinite. Assume it is finite. Then all relations within the postcritical orbit of $F$ over $K$ induce the same kind of relations within the postcritical orbit over any point of~$Y'$. Thus the postcritical orbit over any point of~$Y'$ is finite and there are only finitely many combinatorial possibilities for it. By \cite[Thm.\ 3.3]{Pink2013} or Benedetto-Ingram-Jones-Levy \cite[Cor.~6.3]{BIJL2012} it follows that the number of isomorphism classes of quadratic morphisms arising in all fibers over $Y'$ is finite. Recall that $Y'$ is irreducible of dimension $5$ within~$X$, and that $X$ is finite \'etale over~$U$. Thus the image $V'\subset U$ of $Y'$ is still irreducible of dimension~$5$, and the number of isomorphism classes of quadratic morphisms arising from $F$ in all fibers over $V'$ is finite. 

On the other hand the isomorphism classes of quadratic morphisms over $U$ correspond to the $\GL_{2,k}$-orbits in~$U$ for the action by conjugation on~$F$. Since each orbit has dimension $\le 4$, any irreducible subvariety of dimension $\ge5$ of $U$ must meet infinitely many isomorphism classes. We have therefore obtained a contradiction, and so the postcritical orbit of $F$ over $K$ is infinite, as desired.
\end{Proof}

\begin{Thm}\label{49GarithIncl1}
In the situation above there exists $w\in W$ such that 
$$G^\geom\subset wG(r)w^{-1} \qquad\hbox{and}\qquad 
[G^\arith: G^\arith\cap w G(r)w^{-1}] <\infty.$$
\end{Thm}

\begin{Proof}
Let $k\supset k'\twoheadleftarrow R\into K$ and $F$ be as in Lemma \ref{49pSpecLem}. Then on replacing $k$ by $k'$ the group $G^\geom$ does not change, while $G^\arith$ can only change to a subgroup of finite index. As the desired assertions are invariant under such a modification, we may without loss of generality assume that $k=k'$. 

Let $G_F^\geom \subset G_F^\arith \subset W$ be the geometric and arithmetic monodromy groups of $F$ over~$K$. Then by \cite[Prop.\ \ref{18SpecSpec}]{Pink2013b} we have $G^\geom \subset wG_F^\geom  w^{-1}$ and $G^\arith\subset wG_F^\arith w^{-1}$ for some $w\in W$. On the other hand, let $r'\le r$ be minimal such that $F$ satisfies (\ref{49PostCritRel}) with $r'$ in place of~$r$. Then by Theorems \ref{48GgeomThm} (b) and \ref{48GarithThm1} we have $G_F^\geom=w'G(r')w^{\prime-1}$ and $[G_F^\arith:G_F^\geom]\le 2$. Moreover, by Theorem \ref{47IncluThmAll} we have $G(r')\subset G(r)$. Thus $w'' := w'w \in W$ satisfies $G^\geom\subset w''G(r)w^{\prime\prime-1}$ and $[G^\arith: G^\arith\cap w'' G(r)w^{\prime\prime-1}] \le 2$, and we are done.
\end{Proof}

\medskip
Finally, consider in addition an intermediate field $k\subset k'\subset\bar k$ and a point $x'\in \BP^1(k')\setminus P$, and let $G_{x'} \subset W$ denote the image of the continuous homomorphism $\Gal(\bar k/k')\to W$ describing the action on~$T_{x'}$.

\begin{Thm}\label{49SpecDecomp}
In the situation above there exists $w\in W$ such that 
$$[G^\arith: G^\arith\cap w G(r)w^{-1}] <\infty.$$
\end{Thm}

\begin{Proof}
By \cite[Prop.\ \ref{18SpecDecomp}]{Pink2013b} there exists $w\in W$ such that $G_{x'}\subset wG^\arith w^{-1}$, so the theorem results from the second part of Theorem \ref{49GarithIncl1}.
\end{Proof}

\medskip
This implies a tiny special case of a conjecture of Rafe Jones \cite[Conj. 3.2]{JonesSurvey2013}.

\begin{Ex}\label{49SpecEx}
\rm For any $a\in k\setminus\{0\}$ the quadratic morphism
$$f(x) := \frac{x^2-a}{x^2+a}$$
has the critical points $\infty$ and~$0$, with $f(\infty)=1$ and $f(0)=-1$ and $f^2(\infty) = f^2(0) = \frac{1-a}{1+a}$.
If $a$ is transcendental over the prime field of~$k$, the postcritical orbit of $f$ is infinite by the same argument as in the proof of Lemma \ref{49pSpecLem}. If $k=\BQ$, one checks that the only values of $a$ for which $f$ is conjugate to a quadratic morphism in the list of Manes-Yap \cite[Thm.\ 1.2]{ManesYap} are $a=\pm1$. For these values we have $f^2(\infty) = f^2(0) = 0$ or~$\infty$ and the postcritical orbit of $f$ is finite. For all $a\in \BQ\setminus\{0,\pm1\}$ the postcritical orbit of $f$ is therefore infinite.
\end{Ex}


%
%
%
%
%
%
%


\addcontentsline{toc}{section}{References}




\begin{thebibliography}{99}

\footnotesize

%
%
%
%
%
%
%
%
%
%
%
%
%
%
%
%
%
%
%
%
%
%
%
%
%
%

\bibitem{BIJL2012}
Benedetto, R.L., Ingram, P., Jones, R., Levy, A.:
{\it Critical orbits and attracting cycles in p-adic dynamics}. 
Preprint (Sep 2012) 27p.
{\tt http://arxiv.org/abs/1201.1605}

\bibitem{JonesManes}
Jones, R., Manes. M.:
Galois theory of quadratic rational functions.
to appear in: {\it Commentarii Math. Helvetici}.

\bibitem{JonesSurvey2013}
Jones, R.:
{\it Galois representations from pre-image trees: an arboreal survey.}
In preparation.

\bibitem{ManesYap}
Manes. M., Yap, D.:
{\it A census of quadratic post-critically finite rational maps defined over~$\BQ$.}
Preprint Dec.\ 2012, 18p.\ \ 
{\tt arXiv:1212.1518v1 [math.NT]}

\bibitem{Pink2013}
Pink, R.:
{\it Finiteness and liftability of postcritically finite 
quadratic morphisms in arbitrary characteristic}.
Preprint (version~3, August 2013), 36p.\ \ 
{\tt  arXiv:1305.2841 [math.AG]}

\bibitem{Pink2013b}
Pink, R.:
{\it Profinite iterated monodromy groups
arising from quadratic polynomials}.
Preprint (version~3, September 2013), 85p.\ \ 
{\tt  arXiv:1307.5678 [math.GR]}

\end{thebibliography}
\end{document}